\documentstyle[amscd,amssymb,xypic,12pt]{amsart}
\xyoption{all}
\newcommand{\nc}{\newcommand}

\nc{\exto}[1]{\stackrel{#1}{\longrightarrow}}
\nc{\lan}{\big\langle}
\nc{\ran}{\big\rangle}

\nc{\C}{{\mathsf{k}}}
\nc{\HH}{{\mathbb{H}}}
\nc{\PP}{{\mathbb{P}}}
\nc{\ZZ}{{\mathbb{Z}}}

\nc{\CA}{{\mathcal{A}}}
\nc{\CC}{{\mathcal{C}}}
\nc{\D}{{\mathcal{D}}}
\nc{\CE}{{\mathcal{E}}}
\nc{\CF}{{\mathcal{F}}}
\nc{\CH}{{\mathcal{H}}}
\nc{\CL}{{\mathcal{L}}}
\nc{\CM}{{\mathcal{M}}}
\nc{\CO}{{\mathcal{O}}}
\nc{\CQ}{{\mathcal{Q}}}
\nc{\CU}{{\mathcal{U}}}
\nc{\CV}{{\mathcal{V}}}

\nc{\FE}{{\mathfrak{E}}}
\nc{\FL}{{\mathfrak{L}}}
\nc{\FM}{{\mathfrak{M}}}
\nc{\FS}{{\mathsf{S}}}

\nc{\SK}{{\mathsf{K}}}
\nc{\SO}{{\mathsf{O}}}
\nc{\SX}{{S_X}}
\nc{\SY}{{S_Y}}
\nc{\phipsi}{{q}}

\nc{\UH}{{\mathcal{H}}}

\nc{\TM}{{\widetilde{M}}}
\nc{\TX}{{\widetilde{X}}}
\nc{\TY}{{\widetilde{Y}}}
\nc{\barf}{{\bar{f}}}
\nc{\tf}{{\tilde{f}}}
\nc{\hf}{{\hat{f}}}

\nc{\lotimes}{\mathbin{\mathop{\otimes}\limits^{\mathbb{L}}}}
\nc{\CExt}{\mathop{\mathcal{E}\mathit{xt}}\nolimits}
\nc{\CHom}{\mathop{\mathcal{H}\mathit{om}}\nolimits}
\nc{\RHom}{\mathop{\mathsf{RHom}}\nolimits}
\nc{\RCHom}{\mathop{\mathsf{R}\mathcal{H}\mathit{om}}\nolimits}
\nc{\RG}{\mathop{\mathsf{R\Gamma}}\nolimits}
\nc{\Hom}{\mathop{\mathsf{Hom}}\nolimits}
\nc{\Ext}{\mathop{\mathsf{Ext}}\nolimits}
\nc{\End}{\mathop{\mathsf{End}}\nolimits}
\nc{\Tor}{\mathop{\mathsf{Tor}}\nolimits}
\nc{\Hilb}{\mathop{\mathsf{Hilb}}\nolimits}
\nc{\Spec}{\mathop{\mathsf{Spec}}\nolimits}
\nc{\Pic}{\mathop{\mathsf{Pic}}\nolimits}
\renewcommand{\Im}{\mathop{\mathsf{Im}}\nolimits}
\nc{\Ker}{\mathop{\mathsf{Ker}}\nolimits}
\nc{\Coker}{\mathop{\mathsf{Coker}}\nolimits}
\nc{\codim}{\mathop{\mathsf{codim}}\nolimits}
\nc{\sing}{{\mathsf{sing}}}
\nc{\supp}{{\mathsf{supp}}}
\nc{\rank}{{\mathsf{rank}}}
\nc{\Pf}{{\mathsf{Pf}}}
\nc{\Gr}{{\mathsf{Gr}}}
\nc{\Fl}{{\mathsf{Fl}}}
\nc{\Bl}{{\mathsf{Bl}}}
\nc{\GL}{{\mathsf{GL}}}
\nc{\ev}{{\mathsf{ev}}}
\nc{\id}{{\mathsf{id}}}

\nc{\Cubics}{{{\mathcal{S}}_3}}
\nc{\VFT}{{{\mathcal{S}}_{14}}}
\nc{\VFTE}{{{\mathcal{N}}_{\mathrm{reg,sm}}}}
\nc{\MX}{{\CM_X}}
\nc{\MY}{{\CM_Y}}
\nc{\MYE}{{\CM_{Y,\CE}}}

\theoremstyle{plain}

\newtheorem{theorem}{Theorem}[section]
\newtheorem{lemma}[theorem]{Lemma}
\newtheorem{proposition}[theorem]{Proposition}
\newtheorem{corollary}[theorem]{Corollary}

\theoremstyle{definition}

\newtheorem{definition}[theorem]{Definition}

\theoremstyle{remark}

\newtheorem{remark}[theorem]{Remark}

\newenvironment{proof}{\noindent{\sf Proof:}}{\qed\medskip}

\title{Derived categories of cubic and $V_{14}$ threefolds}
\author{Alexander Kuznetsov}
\address{
Algebra Section, Steklov Mathematical Institute, 
Russian Academy of Sciences, 
 8 Gubkin str., Moscow 119991, Russia
}
\email{akuznet@@mi.ras.ru, sasha@@kuznetsov.mccme.ru}
\date{}

\begin{document}

\maketitle

\centerline{\em In memory of Andrei Nikolaevich Tyurin}

\section{Introduction}

This paper is devoted to the description of several aspects
of a relation of the following two families of Fano threefolds.
The first is the family of cubic threefolds, smooth hypersurfaces
of degree~$3$ in $\PP^4$. The second, is the family
of $V_{14}$ Fano threefolds. It is formed by isomorphism
classes of all smooth complete intersections
$\PP^{9} \cap \Gr(2,6) \subset \PP^{14}$.

The fact that geometry of Fano threefolds from these two families
is related was known for a long time. The history of the question
goes back to Fano himself, who found a birational isomorphism
from a $V_{14}$ threefold to a cubic threefold \cite{Fa,Is}.
Another birational isomorphism was found by Tregub and Takeuchi
\cite{Tr,Ta}.

The paper \cite{IM} has brought a new character into the story,
an instanton bundle on a cubic threefold. An instanton bundle
on a cubic threefold $Y$ is a rank~$2$ stable vector bundle $\CE$
such that $c_1(\CE)=0$ and $H^1(Y,\CE(-1))=0$. 
A topological charge of $\CE$ is defined
as the second Chern class, $c_2(\CE)\in H^4(Y,\ZZ)\cong \ZZ$.
It was shown in~\cite{IM} that for any $V_{14}$ threefold $X$
there exists a unique cubic threefold $Y$ birational to $X$
and that for generic $Y$ the set of $X$ birational to $Y$
is isomorphic to an open subset of the moduli space
$M_0(Y)$ of instanton bundles on $Y$ of topological
charge~$2$.

The goal of the present paper is to show how the above
relation is reflected on the level of the derived categories.
We start however with a more accurate treatment of geometry.
First of all, we remove some genericity coniditions having
been imposed in~\cite{IM} and show that the map $X\mapsto (Y,\CE)$
is actually an isomorphism of moduli stacks. Further, we show
that if $(Y,\CE)$ is the pair, corresponding to $X$,
then we have the following diagram:
$$
\xymatrix{
{\PP_Y(\CE)} \arrow[d]_{p_Y} \arrow[dr]_{\psi} \arrow@{-->}[rr]^{\theta} &&
{\PP_X(\CU)} \arrow[d]^{p_X} \arrow[dl]^{\phi} \\
Y & Q & X
}\eqno{(*)}
$$
where $\CU$ is the restriction of the tautological rank 2 bundle
from the Grassmanian $\Gr(2,6)$ to $X\subset\Gr(2,6)$;
$p_Y$ and $p_X$ are the projectivizations of bundles $\CE$ and $\CU$
over $Y$ and $X$ respectively;
$\psi$ and $\phi$ are small birational contractions onto a
singular quartic hypersurface $Q\subset\PP^5$;
and $\theta = \phi^{-1}\cdot\psi$ is a flop.
The bundle $\CU$ on $X$ is an exceptional bundle.
Thus the above diagram says that the projectivization of the exceptional
bundle on a $V_{14}$ threefold after some natural flop turns into
the projectivization of an instanton bundle on a cubic threefold.

A very similar picture was found in \cite{K} in another situation.
It was shown there that the projectivization of the exceptional
bundle on a $V_{22}$ Fano threefold after a very similar flop turns into
the projectivization of an instanton bundle on the projective space $\PP^3$.
We guess that pictures of this sort should exist for a lot of another
pairs of Fano manifolds and that they are of ultimate importance both 
for the geometry of involved manifolds, and for understanding of Fano 
manifolds in general.

In the second part of the paper we turn our attention to the derived
categories of coherent sheaves on $Y$ and $X$, $\D^b(Y)$ and $\D^b(X)$
respectively. We show that these categories have a similar structure.
First of all, both $\D^b(Y)$ and $\D^b(X)$ contain an exceptional
pair of vector bundles. Explicitly, the pair $(\CO_Y,\CO_Y(1))$
in $\D^b(Y)$, and the pair $(\CO_X,\CU^*)$ in $\D^b(X)$. As usually
in such a situation we obtain semiorthogonal decompositions
$$
\D^b(Y) = \langle \CO_Y, \CO_Y(1), \CA_Y \rangle,\qquad
\D^b(X) = \langle \CO_X, \CU^*, \CA_X \rangle,
$$
where $\CA_Y$ (resp.\ $\CA_X$) is the left orthogonal
to the exceptional pair in $\D^b(Y)$ (resp.\ $\D^b(X)$).
In fact, we use slightly another decomposition of $\D^b(Y)$,
however this change affects only the embedding functor
of $\CA_Y$ into $\D^b(Y)$ and doesn't affect the intrinsic 
structure of $\CA_Y$. Now assume that $Y$ is the cubic threefold 
corresponding to a $V_{14}$ threefold $X$ as above. Then we prove that
the categories $\CA_Y$ and $\CA_X$ are equivalent as triangulated
categories. This is the main result of the paper. The functor, giving
the equivalence is constructed explicitly (see \ref{exp_equ}), 
using diagram~$(*)$.

One of implications of the equivalence is the following. Since all
$V_{14}$ threefolds contained within a fixed birational class
correspond to the same cubic threefold $Y$ it follows that 
the categories $\CA_{X_1}$ and $\CA_{X_2}$ are equivalent 
if $X_1$ and $X_2$ are birational. Thus $\CA_X$ turns into
a birational invariant of~$X$. In fact, we conjecture that
$\CA_X$ allows to distinguish the birational type of $X$,
or equivalently, that $\CA_Y$ allows to distinguish 
the isomorphism class of $Y$. To give some evidence
we construct a family of objects in $\CA_Y$ parameterized
by the Fano surface of lines on $Y$. If one would be able
to describe such a family in intrinsic terms of the category $\CA_Y$ 
(e.g.\ as a moduli space), then it would be possible to reconstruct 
the intermediate Jacobian of $Y$ (as the Albanese variety of the
Fano surface) from $\CA_Y$, and hence, due to the Torelli theorem
\cite{CG,T}, the isomorphism class of~$Y$.

We would like to indicate that the above results can be considered
as a first step to the construction of birational invariants
of algebraic varieties from their derived categories.
We hope this approach might prove useful when dealing 
with the problem of rationality of a cubic fourfold.

\bigskip

The paper is organised as follows. In section~2 we introduce
a definition of the Pfaffian cubic $Y$ and of the theta-bundle $E$,
corresponding to a $V_{14}$ threefold $X$ and state a theorem
on a reconstruction of $X$ from $Y$ and $E$, which is proved
in Appendix~A in a greater generality. After that we introduce
instanton bundles on $Y$ and show that $E$ is a theta-bundle 
iff $E(-1)$ is an instanton of charge~$2$. After that we consider
the projectivizations $\PP_X(\CU)$ and $\PP_Y(E^*)$, construct
their cointractions  
$\phi:\PP_X(\CU) \rightarrow Q \leftarrow \PP_Y(E^*):\psi$
onto a common (singular) quartic hypersurface $Q\subset\PP^5$,
and check that $\theta=\phi^{-1}\circ\psi$ is a flop.
In conclusion we prove some technical results concerning
the fiber product $W = \PP_Y(E^*)\times_Q\PP_X(\CU)$.

We start section~3 with reminding some definitions
and important properties of semiorthogonal decompositions,
mutations, kernel functors, etc. We state also a reformulation
of a result of Bridgeland on flops, which we will need afterwards. 
The remaining part of the section is devoted to the proof of
the main theorem, saying that the categories $\CA_X$ and $\CA_Y$
are equivalent. 

In section~4 we discuss some properties of the category~$\CA_Y$.
First of all, we show that the cube of the Serre functor
of the category $\CA_Y$ is isomorphic to the shift functor,
$\FS_{\CA_Y}^3 \cong [5]$. Moreover, we prove a generalization
of this result for any Fano hypersurface in a projective space.
Also we give examples of two types of objects in $\CA_Y$. 
The first is provided by charge~2 instantons and their shifts, 
and the second is provided by curves on $Y$ with a non-degenerate
theta-characteristics. The particular case of lines on $Y$ gives
a family of objects in $\CA_Y$, parameterized by the Fano surface of $Y$.

In Appendix A we give a general definition of a Pfaffian 
hypersurface and of a theta-bundle and describe some 
of their properties. In Appendix B we give a definition
of instanton bundles on Fano threefolds of index~$2$ and
compute several cohomology groups of their twists.

\subsection*{Notation} 
We assume the base field $\C$ to be an algebraically closed field
of characteristic~$0$. We will use the following notation:
\begin{itemize}\renewcommand{\labelitemi}{-}
\item $V = \C^6$;
\item $A = \C^5$;
\item $f \in \Hom(A,\Lambda^2V^*)$ is an $A$-net of skew-forms on $V$;
\item $X = X_f = \PP(f(A)^\perp) \cap \Gr(2,V)$ 
is a smooth $V_{14}$ Fano threefold;
\item $Y = Y_f = \{ \Pf(f(a))=0 \} \subset \PP(A)$, 
the Pfaffian cubic threefold;
\item $\alpha:Y \to \PP(A)$ is the embedding;
\item $E = E_f$ is the theta-bundle on $Y$;
\item $\CE$ is an instanton of charge~$2$ on $Y$, $\CE = E(-1)$;
\item $\CU$ is a restriction of the tautological vector bundle 
from $\Gr(2,V)$ to $X$;
\item $p_X:\PP_X(\CU) \to X$ is the projectivization of $\CU$ on $X$;
\item $p_Y:\PP_Y(E^*) \to Y$ is the projectivization of $E^*$ on $Y$;
\item $\phi:\PP_X(\CU) \to \PP(V)$ is the map, 
induced by embedding $\PP_X(\CU) \subset \Fl(1,2;V)$;
\item $\psi:\PP_Y(E^*) \to \PP(V)$ is the map, 
induced by embedding $\PP_Y(E^*) \subset \Fl(1,2;V)$;
\item $Q = \phi(\PP_X(\CU)) = \psi(\PP_Y(E^*)) \subset \PP(V)$ 
is a quartic hypersurface;
\item $C = \sing(Q)$ is a curve, $\deg C = 25$, $p_a(C)=26$;
\item $\SX \subset \PP_X(\CU)$ 
is a ruled surface, contracted by $\phi$ to $C$;
\item $\SY \subset \PP_Y(E^*)$ 
is a ruled surface, contracted by $\psi$ to $C$;
\item $\theta: \PP_Y(E^*) \to \PP_X(\CU)$ is the flop in $\SY$, 
$\theta = \phi^{-1}\circ\psi$;
\item $W = \PP_Y(E^*) \times_Q \PP_X(\CU)$ is the fiber product;
\item $\eta:W \to \PP_Y(E^*)$, $\xi:W \to \PP_X(\CU)$ and $\phipsi:W \to Q$
are the projections;
\item $i:W \to \PP_Y(E^*) \times \PP_X(\CU)$, $j:W \to Y\times X$
and $\lambda: W \to \PP(A)\times\PP_X(\CU)$ are the embeddings.
\end{itemize}

\subsection*{Acknowledgements}
I am grateful to Dmitry Orlov and Alexei Bondal for useful discussions.
I was partially supported by RFFI grants 02-01-00468 and 02-01-01041 
and INTAS-OPEN-2000-269. The research described in this work 
was made possible in part by CRDF Award No. RM1-2406-MO-02.
A part of this work was accomplished during my visit at
the Universite Paul Sabatier (Toulouse) and Insitute de Mathematique
de Luminy (Marseille), which was organised by the National Scientific 
Research Center of France and by the Independent University 
of Moscow via the ``Jumelage Mathematique'' program.

Finally, I would like to express my sincerest grattitude
to Andrei Nikolaevich Tyurin whose ideas always were a source
of inspiration and whose work was an object of admiration for me.

\section{Geometry}

Consider a five-dimensional vector space $A=\C^5$,
a six-dimensional vector space $V=\C^6$, and a
linear map $f:A\to\Lambda^2V^*$. Such map is called
an {\em $A$-net of skew-forms on $V$}.

\subsection*{Pfaffian cubic and theta-bundle}

For any such $f$ let $f(A)^\perp\subset\Lambda^2V$ denote
the annihilator of $f(A)\subset\Lambda^2V^*$. Denote also
$X = X_f = \PP(f(A)^\perp) \cap \Gr(2,V) \subset \PP(\Lambda^2V)$.
When $f$ is generic $X$ is a smooth Fano threefold of index 1
with $\Pic X=\ZZ$ and of genus~$8$. Such threefolds are known as
{\em $V_{14}$ Fano threefolds}\/ \cite{Is1,IP}. Moreover, any $V_{14}$
threefold can be realized as $X_f$ for some $f$ \cite{Mu}.

An $A$-net $f$ is called {\em regular}\/ if $\rank f(a) \ge 4$
%the rank of $f(a)$ is not less then~$4$ 
for any $0\ne a\in A$.

\begin{lemma}\label{reg}
If $X = X_f \subset \Gr(2,V)$ is a smooth $V_{14}$ threefold
then the $A$-net $f$ is regular.
\end{lemma}
\begin{proof}
Assume that the $A$-net $f:A\to\Lambda^2V^*$ isn't regular.
Then the rank of a skew-form $f(a)\in\Lambda^2 V^*$ is less
or equal than~$2$ for some $0\ne a\in A$. Let $K_a = \Ker f(a) \subset V$
be the kernel of this form. Then $\dim K_a\ge 4$ and the Grassmannian
$\Gr(2,K_a)\subset\Gr(2,V)$ has nonempty intersection with $X$, because
$X\cap\Gr(2,K_a)$ is a plane section of $\Gr(2,K_a)$ of codimension $\le 4$,
and $\dim\Gr(2,K_a)\ge 4$. But it is easy to check that any point in
$X\cap\Gr(2,K_a)$ is singular in $X$ (see the proof of proposition~\ref{sm}).
\end{proof}

Any $A$-net $f$ can be considered as an element of
$\Hom(V\otimes\CO_{\PP(A)}(-1),V^*\otimes\CO_{\PP(A)})$,
the space of homomorphisms of coherent sheaves on $\PP(A)$.
If $f$ is regular then this homomorphism is injective, and
its cokernel $E=E_f$ is a sheaf supported on a cubic hypersurface
$Y=Y_f \subset \PP(A)$ with equation $\Pf(f(a)) = 0$ (where $\Pf$ 
stands for the Pfaffian of a skew-form), the {\em Pfaffian cubic}\/ of $f$.
Thus we have an exact sequence of coherent sheaves on $\PP(A)$:
\begin{equation}\label{es}
0 \to V\otimes\CO_{\PP(A)}(-1) \exto{f} 
V^*\otimes\CO_{\PP(A)} \to \alpha_*E \to 0,
\end{equation}
where $\alpha:Y\to\PP(A)$ is the embedding. We call $E_f$
{\em the theta-bundle}\/ of the $A$-net $f$ (see Appendix~A).
The map $V^*\otimes\CO_{\PP(A)} \to \alpha_*E$ induces an isomorphism
$\gamma_f:V^* = H^0(\PP(A),V^*\otimes\CO_{\PP(A)}) \to 
H^0(\PP(A),\alpha_*E) = H^0(Y_f,E_f)$.

\begin{theorem}\label{v14}
Associating to an $A$-net $f$ the triple $(Y_f,E_f,\gamma_f)$ 
gives a $\GL(A)\times\GL(V)$-equivariant
isomorphism between
\begin{itemize}
\item the subset of $\PP(A^*\otimes\Lambda^2V^*)$ formed by all
regular $A$-nets of skew-forms on $V$, and
\item the set of triples $(Y,E,\gamma)$, where
$Y$ is a cubic hypersurface in $\PP(A)$,
$E$ is a rank~$2$ locally free sheaf on $Y$,
and $\gamma$ is an isomorphism $V^*\to H^0(Y,E)$,
such that
\begin{equation}\label{chern_and_h}
c_1(E) = 2[h],\quad c_2(E) = 5[l],\quad\text{and}\quad
H^\cdot(Y,E(t)) = 0\quad\text{for $-3\le t\le -1$},
\end{equation}
where $[h]\in H^2(Y,\ZZ)$ and $[l]\in H^4(Y,\ZZ)$ are
the classes of a hyperplane section and of a line
respectively.
\end{itemize}
Further, the theta-bundle $E_f$ of
a regular $A$-net is generated by global sections,
$H^0(Y_f,E_f)=V^*$, and induces an embedding
$\kappa:Y_f\to\Gr(2,V)$. Finally,
$\sing(X_f) = \sing(Y_f) = X_f\cap Y_f \subset \Gr(2,V)$.
In particular, $Y_f$ is smooth iff $X_f$ is smooth.
\end{theorem}

The major part of this theorem is proved in \cite{MT,IM,Beau,Dr}
in more or less the same generality. Only the last statement
seems to be new. We give a complete proof in Appendix A.

\begin{remark}\label{ffromx}
It is easy to check that $H^0(X_f,\CO_{X_f}(1)) \cong \Lambda^2V^*/f(A)$.
It follows that the $A$-net $f$ can be reconstructed from $X_f$ up to
the action of $\GL(A)\times\GL(V)$, the action of $\GL(V)$ corresponds 
to a choice of embedding $X_f \to \Gr(2,V)$, and the action of $\GL(A)$
corresponds to a choice of isomorphism 
$A\to\Ker(\Lambda^2V^* = H^0(\Gr(2,V),\CO_{\Gr(2,V)}(1)) \to 
H^0(X_f,\CO_{X_f}(1)))$.
\end{remark}

If $X$ is a smooth $V_{14}$ Fano threefold and $f$ is an $A$-net 
of skew-forms on $V$, such that $X\cong X_f$, then by remark~\ref{ffromx}
and theorem~\ref{v14}, the pair $(Y_f,E_f)$ is determined by $X$
up to an isomorphism. We will say that $Y_f$ is the Pfaffian cubic
of $X$ and $E_f$ is the corresponding theta-bundle.

\subsection*{Instantons}

\begin{definition}
A sheaf $\CE$ on a cubic threefold $Y\subset\PP^4$ is
an {\em instanton bundle}\/ if $\CE$ is locally free
of rank~$2$, stable (with respect to $\CO_Y(1)$) and
$c_1(\CE) = 0$, $H^1(Y,\CE(-1))=0$.
The {\em topological charge}\/ of an instanton $\CE$ is an integer $k$,
such that $c_2(\CE)=k[l]$, where $[l]\in H^4(Y,\ZZ)$ is the class of a line.
\end{definition}

This definition is a straightforward generalization of the definition
of (mathematical) instanton vector bundle on $\PP^3$ \cite{OSS} and admits
further generalization to any Fano threefold of index~$2$.
We introduce such definition and deduce simplest implications
in Appendix~B. It is shown, in particular, that the smallest possible
charge for the instantons on $Y$ is~$2$, and

\begin{proposition}\label{charge2}
If $\CE$ is an instanton vector bundle of charge~$2$ on $Y$ then
$$
H^p(Y,\CE(t)) = \begin{cases}
\C^6, & \text{for $(p,t)=(0,1)$ and $(p,t)=(3,-3)$}\\
0  , & \text{for other $(p,t)$ with $-3\le t\le 1$}
\end{cases}
$$
\end{proposition}

Consider, following \cite{MT} the Gieseker--Maruyama moduli space
$M_Y(2;0,2)$ of semistable (with respect to $\CO_Y(1)$) rank~2
torsion free sheaves on $Y$ with Chern classes $c_1=0$ and $c_2=2[l]$
and its Zariski open subset
\begin{multline}\label{m0}
M_0(Y) =\{[\CE ]\in M_Y(2;0,2) \; |\; (i)\; \CE\;\mbox{is
{\em stable}\/ and {\em locally free}};\\ (ii) \; H^1(Y,\CE (-1))=
H^1(Y,\CE (1))=H^2(Y,\CE (1))=H^2(Y,\CE\otimes\CE )=0
\} .
\end{multline}

\begin{proposition}\label{pm0}
The following conditions are equivalent:
$$
\begin{array}{rl}
(i)\qquad   & \text{$\CE$ is an instanton bundle of charge $2$;}\smallskip\\
(ii)\qquad & \text{$\CE(1)$ satisfies conditions $(\ref{chern_and_h})$;}\smallskip\\
(iii)\qquad & \text{$\CE(1)$ is a theta-bundle;}\smallskip\\
(iv)\qquad  & \text{$[\CE]\in M_0(Y)$.}
\end{array}
$$
\end{proposition}
\begin{proof}
The implication $(i) \Rightarrow (ii)$ easily follows from
proposition~\ref{charge2}, $(ii) \Rightarrow (iii)$ is given by
theorem~\ref{v14}, $(iv) \Rightarrow (i)$ is trivial. 
Thus it remains to check the implication $(iii) \Rightarrow (iv)$.

Assume that $E$ is a theta-bundle of $f$ and denote $\CE=E(-1)$.
It follows from theorem~\ref{v14} that it suffices to check
that $H^2(Y,\CE\otimes\CE)=0$. Restricting (\ref{es}) to $Y$
and taking into account isomorphism
$L^1\alpha^*\alpha_*E \cong 
E\otimes L^1\alpha^*\alpha_*\CO_Y \cong 
E\otimes\CO_Y(-3)$
we get the following exact sequence
\begin{equation*}\label{esy}
0 \to E(-3) \to V\otimes\CO_Y(-1) \to V^*\otimes\CO_Y \to E \to 0.
\end{equation*}
Applying $\Hom(E,-)$ and taking into account isomorphisms
$$
\begin{array}{llllll}
\Ext^p(E,\CO_Y)     &\cong& H^p(Y,E^*)     &\cong& H^p(Y,E(-2)) &= 0,\\
\Ext^p(E,\CO_Y(-1)) &\cong& H^p(Y,E^*(-1)) &\cong& H^p(Y,E(-3)) &= 0,
\end{array}
$$
(we used here an isomorphism $\det E^*\cong\CO_Y(-2)$ and 
properties (\ref{chern_and_h})) we obtain isomorphisms
\begin{equation}\label{extee}
\Ext^p(E,E) \cong \Ext^{p+2}(E,E(-3)).
\end{equation}
In particular, $\Ext^2(E,E)=0$, but
$\Ext^2(E,E)\cong H^2(Y,E^*\otimes E) \cong H^2(Y,\CE\otimes\CE)$.
\end{proof}

\begin{remark}
Note, that by the way we have proved that conditions
$$
H^1(Y,\CE(1)) = H^2(Y,\CE(1)) = H^2(Y,\CE\otimes\CE) = 0
$$
in the definition of the set $M_0(Y)$ in~\cite{MT} are redundant.
\end{remark}

\begin{remark}\label{es_fib}
The embedding $E^*(-1) \cong E(-3) \to V\otimes\CO_Y(-1)$
identifies the fiber of $E^*(-1)$ over a point $a\in Y$
with $\Ker f(a)\subset V$.
\end{remark}

\subsection*{Moduli stacks}

Let $\VFTE$ denote the open subset of $\PP(A^*\otimes\Lambda^2 V^*)$,
consisting of regular $A$-nets $f$, such that $X_f$ (or, equivalently,
$Y_f$) is smooth. Let $\MX$ denote the moduli stack of $V_{14}$
threefolds. Let $\MY$ denote the moduli stack of cubic threefolds.
Let $\MYE$ denote the moduli stack of pairs $(Y,\CE)$, where $Y$ is
a cubic threefold and $\CE$ is an instanton of charge~$2$ on $Y$.

\begin{theorem}\label{modulistacks}
We have an isomorphism of stacks
$\MX \cong \MYE \cong \VFTE/\!/\GL(A)\times\GL(V)$.
In particular, the fiber of the projection $\MX\to\MY$ 
is isomorphic to $M_0(Y)$.
\end{theorem}
\begin{proof}
It suffices to repeat the arguments of theorem~\ref{v14} 
and remark~\ref{ffromx}
% and of proposition~\ref{pm0} 
in relative situation.
% and to apply the following.
\end{proof}

\subsection*{Further properties of theta-bundles}

\begin{lemma}\label{temy}
If $E$ is a theta-bundle then $H^p(Y,S^2 E(-1))=
\begin{cases}\C, & \text{if $p=1$}\\0, & \text{otherwise}\end{cases}$
\end{lemma}
\begin{proof}
It follows from (\ref{extee}) that
$$
\Hom(E,E(-3)) = \Ext^1(E,E(-3)) = 0
\qquad\text{and}\qquad
\Ext^2(E,E(-3)) = \Hom(E,E) = \C,
$$ 
because $E$ is stable. 
Applying the Riemann-Roch we deduce that $\Ext^3(E,E(-3)) = \C^5$. 
Further, taking into account an isomorphism $E^*(1)\cong E(-1)$
and applying the Serre duality on $Y$
we obtain
$$
H^p(Y,E\otimes E(-1)) \cong
%H^p(Y,E\otimes E^*(1)) \cong
\Ext^p(E^*(1),E) \cong
\Ext^p(E(-1),E) \cong
\Ext^{3-p}(E,E(-3))^* \cong
\begin{cases}
\C^5, & \text{if $p=0$}\\
\C  , & \text{if $p=1$}\\
0   , & \text{otherwise}
\end{cases}
$$
But $E\otimes E(-1)\cong \Lambda^2 E(-1) \oplus S^2 E(-1) \cong
\CO_Y(1) \oplus S^2 E(-1)$.
So, it remains to note that
$H^0(Y,\CO_Y(1))=H^0(\PP^4,\CO_{\PP^4}(1)) \cong \C^5$, 
$H^{>0}(Y,\CO_Y(1))=0$,
and  lemma follows.
\end{proof}

\subsection*{$\PP^1$-bundle over $X$}

Let $X\subset\Gr(2,V)$ be a smooth $V_{14}$ Fano threefold
and let $f:A\to\Lambda^2V^*$ be the corresponding $A$-net.
Let $\CU$ denote the restriction to $X$ of the tautological rank~2 subbundle
on $\Gr(2,V)$. Then the projectivization $p_X:\PP_X(\CU)\to X$ is embedded
into the partial flag variety $\Fl(1,2;V)$. Let $\phi:\PP_X(\CU)\to\PP(V)$
denote the restriction of the canonical projection $\Fl(1,2;V)\to\PP(V)$.

\begin{proposition}\label{phi}
$(i)$ The image $Q=\phi(\PP_X(\CU))\subset\PP(V)$ is a quartic hypersurface,
singular along a curve $C\subset Q$, $\deg C=25$, $p_a(C)=26$.
$(ii)$ The map $\phi:\PP_X(\CU) - \phi^{-1}(C)\to Q - C$ is an isomorphism,
while $\phi:\phi^{-1}(C)\to C$ is a $\PP^1$-bundle.
$(iii)$ For any point $c\in C$ the curve $L_c=p_X(\phi^{-1}(c))\subset X$
is a line on $X$. On the other hand, if $L$ is a line on $X\subset\Gr(2,V)$,
then $p_X^{-1}(L)$ is a Hirzebruch surface $F_1$, and its exceptional
section $\tilde{L}$ coincides with $\phi^{-1}(c)$ for some $c\in C$.
\end{proposition}
\begin{proof}
$(i)$ It is clear that the image $Q=\phi(\PP_X(\CU))\subset\PP(V)$ is
just the set of all points $v\in\PP(V)$ which are contained in a
2-dimensional subspace $U\subset V$ isotropic with respect to all
skew-forms in the $A$-net $f$. Thus $v\in Q$ if and only if the map
$$
f_v:A \to V^*,\qquad a\mapsto f(a)(v,-)
$$
has image of codimension $\ge 2$. In the other words, $Q$ is the
determinantal $\{\rank(f)\le 4\}\subset\PP(V)$, where $f$
is considered as a homomorphism of coherent sheaves on $\PP(V)$
$$
A\otimes\CO_{\PP(V)}(-1) \exto{f} V^*\otimes\CO_{\PP(V)}.
$$
Note that since $f(a)$ is a skew-form we have $f(a)(v,v)=0$ for all
$a\in A$, hence the image of $f$ lies in
the annihilator $v^\perp\subset V^*$. In the other words,
the above homomorphism of sheaves factors as
$$
A\otimes\CO_{\PP(V)}(-1) \exto{f'}
\Omega_{\PP(V)}(1)
\subset V^*\otimes\CO_{\PP(V)}
$$
Note that $\rank(A\otimes\CO_{\PP(V)}(-1)) = \rank(\Omega_{\PP(V)}(1)) = 5$,
hence $Q$ is the zero locus of
$$
\det f'\in\Hom(\det(A\otimes\CO_{\PP(V)}(-1)),\det(\Omega_{\PP(V)}(1))) =
\Hom(\CO_{\PP(V)}(-5),\CO_{\PP(V)}(-1)).
% = \Gamma(\PP(V),\CO_{\PP(V)}(4)).
$$
Thus $Q$ is a quartic hypersurface in $\PP(V)$.

By the general properties of determinantals the singular locus
$C=\sing(Q)$ is the determiantal $\{\rank(f')\le 3\}\subset\PP(V)$.
It will be shown in $(iii)$ below that $C$ parameterizes lines on $X$,
hence $C$ is $1$-dimensional \cite{Is1,IP}. Thus $C$ is of expected dimension 
and the standard methods can be applied to compute $\deg(C)=25$, $p_a(C)=26$.

$(ii)$ It is clear that for a point $v\in Q$ the fiber $\phi^{-1}(v)$
is isomorphic to the set of all $2$-dimensional subspaces $U\subset V$,
such that $v\in U$ and $U\subset(\Im f_v)^\perp$.
But $\dim (\Im f_v)^\perp = 2$ for $v\in Q - C$, hence $\phi$ is
an isomorphism over the complement of $C$. On the other hand,
for any point $v\in C$ we have
$(\Im f_v)^\perp \cong \C^3$ and 
$\phi^{-1}(v) \cong \PP((\Im f_v)^\perp/\C v) \cong \PP^1$.
%
%$$
%\phi^{-1}(v) = \{ \C^2 = U\subset V\ |\
%\langle v\rangle \subset U\subset (\Im f_v)^\perp \cong \C^3
%\subset V \} \cong \PP^1.
%$$

$(iii)$ The arguments in $(ii)$ show that $L_c:=p_X(\phi^{-1}(c))$ 
is a line on $X\subset\Gr(2,V)$ for any $c\in C$. On the other hand,
if $L$ is a line on $X\subset\Gr(2,V)$, then 
$\CU_{|L}\cong \CO_L\oplus\CO_L(-1)$. Thus
$$
p_X^{-1}(L) = \PP_L(\CU_{|L}) = \PP_L(\CO_L\oplus\CO_L(-1)) \cong F_1.
$$
It is clear that the map $\phi$ contracts the exceptional 
section $\tilde{L}$ of $p_X^{-1}(L)$, hence $\tilde{L}=\phi^{-1}(c)$ 
for some $c\in C$ and $L = L_c$.
\end{proof}

\begin{remark}\label{rge3}
It is clear that for any $0\ne v\in V$ we have $\PP(\Ker f_v)\subset Y_f$.
On the other hand $Y_f$ is a smooth cubic in $\PP(A)$ by theorem~\ref{v14},
hence it cannot contain a $\PP^2$. This means that $\dim\Ker f_v\le 2$
and $\rank(f_v)\ge 3$ for any $0\ne v\in V$.
\end{remark}

\begin{remark}
Using description of $C$ as a determinantal one can show
that $\CO_{\PP(V)}(1)_{|C}$ is a degenerate even theta-characteristic
on $C$ with $\dim H^0(C,\CO_{\PP(V)}(1)_{|C}) = 6$.
\end{remark}

\begin{corollary}\label{lines_x}
The curve $C$ parameterizes lines on $X$.
\end{corollary}

\subsection*{$\PP^1$-bundle over $Y$}

Now let $Y$ be the Pfaffian cubic of $X$ and
let $E$ be the theta-bundle of $X$. By theorem~\ref{v14} 
the bundle $E$ induces an embedding $\kappa:Y\to\Gr(2,V)$.
Then we obtain an embedding of the projectivization $p_Y:\PP_Y(E^*)\to Y$ 
into the partial flag variety $\Fl(1,2;V)$. Let $\psi:\PP_Y(E^*)\to\PP(V)$
denote the restriction of the canonical projection $\Fl(1,2;V)\to\PP(V)$.

\begin{proposition}\label{psi}
$(i)$ We have $\psi(\PP_Y(E^*)) = Q$.
$(ii)$ The map $\psi:\PP_Y(E^*) - \psi^{-1}(C)\to Q - C$ is an isomorphism,
while $\psi:\psi^{-1}(C)\to C$ is a $\PP^1$-bundle.
$(iii)$ For any point $c\in C$ the curve $M_c=p_Y(\psi^{-1}(c))\subset Y$
is a line on $Y$ such that $E^*_{|M_c}\cong \CO_{M_c}\oplus\CO_{M_c}(-2)$.
On the other hand, if $M$ is a line on $Y$ such that
$E^*_{|M}\cong \CO_{M}\oplus\CO_{M}(-2)$, then $p_Y^{-1}(M)$ 
is a Hirzebruch surface $F_2$, and its exceptional
section $\tilde{M}$ coincides with $\psi^{-1}(c)$ for some $c\in C$.
\end{proposition}
\begin{proof}
$(i)$ The fiber of $E^*$ over a point $a\in Y$ is the kernel
of the skew-form $f(a)\in\Lambda^2V^*$. Hence $v\in\psi(\PP_Y(E^*))$
iff $f(a)(v,-)=0$, that is iff $f_v(a)=0$ for some $0\ne a\in A$.
Thus $\psi(\PP_Y(E^*)) = Q$.

$(ii)$ Note that the fiber of $\psi$ over $v\in\PP(V)$
coincides with $\PP(\Ker f_v)\subset\PP(A)$.
For $v\in Q - C$ we have $\rank(f_v)=4$, hence $\dim\Ker f_v=1$.
Thus $\psi$ is an isomorphism over $Q - C$.
On the other hand, for $v\in C$ we have $\rank(f_v)=3$,
hence $\dim\Ker f_v=2$ and $\psi^{-1}(v)\cong\PP^1$.

$(iii)$ The arguments in $(ii)$ show that $M_c=p_Y(\psi^{-1}(c))$
is a line on the cubic $Y\subset\PP(A)$. Note that
$\det(E^*_{|M_c}) = \det(E^*)_{|M_c} = \CO_Y(-2)_{|M_c} = \CO_{M_c}(-2)$,
and $c\in V$ gives a nonvanishing section of 
$E^*_{|M_c}\subset V\otimes\CO_{M_c}$, hence
$E^*_{|M_c}\cong\CO_{M_c}\oplus\CO_{M_c}(-2)$. 
On the other hand, if $M$ is a line on $Y$ such that
$E^*_{|M}\cong \CO_{M}\oplus\CO_{M}(-2)$, then 
$$
p_Y^{-1}(M) = \PP_{M_c}(E^*_{|M}) =
\PP_{M}(\CO_{M}\oplus\CO_{M}(-2)) \cong F_2.
$$
It is clear that the map $\psi$ contracts the exceptional 
section $\tilde{M}$ of $p_Y^{-1}(M)$, hence $\tilde{M}=\psi^{-1}(c)$ 
for some $c\in C$ and $M = M_c$.
\end{proof}

\begin{corollary}\label{jlines_y}
The curve $C$ parameterizes jumping lines of the instanton $\CE=E(-1)$ on $Y$.
\end{corollary}

\subsection*{The flop}

Denote $\SX=\phi^{-1}(C)\subset\PP_X(\CU)$ and
$\SY=\psi^{-1}(C)\subset\PP_Y(E^*)$. Thus $\SX$ and $\SY$ are ruled
surfaces over the curve $C$. It is proved in propositions~\ref{phi}
and~\ref{psi} that $\phi$ contracts $\SX$ onto~$C$ and $\psi$
contracts $\SY$ onto~$C$. Hence the rational map
$\theta=\phi^{-1}\circ\psi : \PP_Y(E^*) \to \PP_X(\CU)$
is a birational isomorphism.

\begin{theorem}\label{flop}
The map $\theta$ is a flop in the surface $\SY$.
The map $\theta^{-1}$ is a flop in the surface $\SX$.
\end{theorem}
\begin{proof}
Since $\phi$ and $\psi$ are small contractions by
propositions~\ref{phi} and \ref{psi}, it remains to check
that the canonical classes of $\PP_Y(E^*)$ and $\PP_X(\CU)$
are pull-backs from $Q$. But it is easy to see that the canonical classes
equal $\psi^*\CO_Q(-2)$ and $\phi^*\CO_Q(-2)$ respectively. Indeed,
$$
\arraycolsep = 2pt
\begin{array}{llll}
\omega_{\PP_Y(E^*)} & =
p_Y^*\omega_Y \otimes \omega_{\PP_Y(E^*)/Y} & =
p_Y^*\CO_Y(-2) \otimes \big(\psi^*\CO_Q(-2) \otimes p_Y^*\det E\big) & \cong
\psi^*\CO_Q(-2),\\
\omega_{\PP_X(\CU)} & =
p_X^*\omega_X \otimes \omega_{\PP_X(\CU)/X} & =
p_X^*\CO_X(-1) \otimes \big(\phi^*\CO_Q(-2) \otimes p_X^*\det \CU^*\big) & \cong
\phi^*\CO_Q(-2).
\end{array}
$$
since $\psi^*\CO_Q(1)$ and $\phi^*\CO_Q(1)$ are the Grothendieck
relatively ample line bundles on $\PP_Y(E^*)$ and $\PP_X(\CU)$ respectively
by definition of $\psi$ and $\phi$.
\end{proof}

Summarizing, we get the following.

\begin{theorem}\label{main}
Let $X$ be a smooth $V_{14}$ Fano threefold.
Let $Y$ be its Pfaffian cubic and let $E$ be the theta-bundle of $X$ on $Y$.
Then we have the following diagram
$$
\xymatrix{
{\PP_Y(E^*)} \arrow[dd]_{p_Y} \arrow[ddrr]^{\psi} 
\arrow@/^4ex/@{-->}[rrrr]^{\theta} &
{\ \SY\ } \arrow@{_{(}->}[l] \arrow[dr] &&
{\ \SX\ } \arrow@{^{(}->}[r] \arrow[dl] &
{\PP_X(\CU)} \arrow[dd]^{p_X} \arrow[ddll]_{\phi} \\
&& {\mathstrut C} \arrow@{^{(}->}[d] \\
Y && Q && X
}
$$
where
\begin{itemize}
\item $Q$ is a quartic hypersurface in $\PP(V)$,
singular along a curve $C$;
\item $\SX\subset\PP_X(\CU)$ and $\SY\subset\PP_Y(E^*)$
are ruled surfaces over the curve $C$, ruled by exceptional
sections over lines on $X$ and by exceptional sections over
jumping lines on $Y$ respectively;
\item $\phi$ and $\psi$ contract ruled surfaces $\SX$ and $\SY$ onto $C$
and bijective elsewhere;
\item $\theta = \phi^{-1}\cdot\psi$ is a flop in $\SY$.
\end{itemize}
\end{theorem}

\begin{remark}[\cite{IM}]\label{xy_bir}
If $H$ is a hyperplane in $\PP(V)$ then it is easy to see that
$p_X\circ\phi^{-1}:Q\cap H \to X$ and $p_Y\circ\psi^{-1}:Q\cap H \to Y$
are birational isomorphisms. In particular, the Pfaffian cubic $Y$
of a smooth $V_{14}$ Fano threefold $X$ is birational to $X$.
Moreover, the Torelly theorem \cite{CG,T} implies that 
cubic threefolds $Y_1$ and $Y_2$ are birational if and only if
they are isomorphic. It follows that the fibers of the map
of the moduli stacks $\MX\to\MY$ are birational classes 
of $V_{14}$ threefolds.
\end{remark}

\subsection*{The fiber product}

Consider the fiber product $W=\PP_Y(E^*)\times_Q\PP_X(\CU)$ 
and denote the embedding $W \to \PP_Y(E^*)\times\PP_X(\CU)$ by~$i$.
Let $\xi:W \to \PP_X(\CU)$, $\eta:W \to \PP_Y(E^*)$, $\phipsi: W\to Q$
denote the projections.
Put $j=(p_Y\times p_X)\cdot i:W \to Y\times X$ and
$\lambda = (\alpha p_Y\times\id)\cdot i : W \to \PP(A)\times\PP_X(\CU)$.
$$
\xymatrix{
W \arrow[r]^\xi \arrow[d]_\eta \arrow[dr]^\phipsi &
\PP_X(\CU) \arrow[d]^\phi \\
\PP_Y(E^*) \arrow[r]^\psi & 
Q
}
\qquad
\xymatrix{
&& W \arrow[lld]_-{j} \arrow[d]^{i} \arrow[rrd]^-{\lambda} \\
Y\times X &&
\PP_Y(E^*)\times\PP_X(\CU) 
\arrow[ll]_-{p_Y\times p_X} \arrow[rr]^-{\alpha p_Y\times \id} &&
\PP(A)\times\PP_X(\CU)
}
$$

\begin{proposition}\label{fp}
$(i)$ 
$j$ is a closed embedding and we have the following
exact sequence on $Y\times X$:
\begin{equation}\label{jw}
0 \to
E^*(-1) \boxtimes \CO_X \to
\CO_Y(-1)\boxtimes V/\CU \exto{}
\CO_Y\boxtimes \CU^* \exto{}
j_*\phipsi^*\CO_Q(1) \to 0,
\end{equation}
$(ii)$ 
$\lambda$ is a closed embedding and we have the following
exact sequence on $\PP(A)\times\PP_X(\CU)$:
\begin{equation}\label{jw1}
\begin{array}{r}
0                                                                     \to
\CO_{\PP(A)}(-4)\boxtimes\Lambda^4(\phi^*\CO_Q(-1)\otimes p_X^*V/\CU) \to
\CO_{\PP(A)}(-3)\boxtimes\Lambda^3(\phi^*\CO_Q(-1)\otimes p_X^*V/\CU) \to \\ \to
\CO_{\PP(A)}(-2)\boxtimes\Lambda^2(\phi^*\CO_Q(-1)\otimes p_X^*V/\CU) \to
\CO_{\PP(A)}(-1)\boxtimes\Lambda^1(\phi^*\CO_Q(-1)\otimes p_X^*V/\CU) \to \\ \to
\CO                                                                   \to
\lambda_*\CO_W                                                        \to
0.
\end{array}
\end{equation}
\end{proposition}
\begin{proof}
$(i)$
By definition of $X$ the composition
$\CU \to V\otimes\CO_X \exto{f(a)} V^*\otimes\CO_X \to \CU^*$
vanishes for any $a\in A$. Therefore, $f$ induces a morphism 
of vector bundles
$\CO_Y(-1)\boxtimes V/\CU \exto{\tf} \CO_Y\boxtimes \CU^*$.
Moreover, in the following commutative diagram
$$
\xymatrix{
\Ker f(a) \arrow[r] \arrow[dr] &
V\otimes\CO_X \arrow[r]^{f(a)} \arrow[d] &
V^*\otimes\CO_X \arrow[d] \arrow[r] &
\Coker f(a) \arrow[d] \\
& V/\CU \arrow[r]^{\tf_a} & \CU^* \arrow[r] & \Coker \tf_a
}
$$
the upper row is a complex, hence the sequence
\begin{equation}\label{sq3}
E^*(-1) \boxtimes \CO_X \to
\CO_Y(-1)\boxtimes V/\CU \exto{\tf}
\CO_Y\boxtimes \CU^*,
\end{equation}
in which the first morphism is the canonical embedding
$E^*(-1) \boxtimes \CO_X \to \CO_Y(-1)\boxtimes V/\CU$
(corresponding to the diagonal arrow in the diagram, 
see remark~\ref{es_fib}), is a complex.

For any point $(a,U)\in Y\times X \subset \PP(A)\times\Gr(2,V)$
the composition $\Ker f(a) \to V \to V/U$ is an embedding unless
$\Ker f(a)\cap U\ne 0$. Similarly, the map $V/U \exto{\tf_a} U^*$
is a surjection unless $\Ker f(a)\cap U\ne 0$.
Any $0\ne v\in\Ker f(a)\cap U$ specifies a point 
$(a,v)\in\PP_Y(E^*)\subset\PP(A)\times\PP(V)$
and a point $(U,v)\in\PP_X(\CU)\subset\Gr(2,V)\times\PP(V)$ 
such that $\psi(a,v)=v=\phi(U,v)\in\PP(V)$.
This means that the degeneration sets of both morphisms
of~(\ref{sq3}) coincide with $j(W)$.

Let $W'$ denote the degeneration subscheme of the morphism $\tf$
on $Y\times X$. We already have shown that $j$ is a set-theoretical
bijection $W\to W'$. Let us show that $j$ is a scheme-theoretical isomorphism.

Indeed, the pullback of the morphism $\tf$ via $j$ is
$\eta^*p_Y^*\CO_Y(-1) \otimes \xi^*p_X^*V/\CU \to \xi^*p_X^*\CU^*$.
It is clear that its composition with the surjection
$(\xi^*p_X^*\CU^* \to \phipsi^*\CO_Q(1)) = 
\xi^*(p_X^*\CU^* \to \phi^*\CO_Q(1))$
vanishes, hence $j:W\to Y\times X$ factors through
the subscheme $W'\subset Y\times X$.

On the other hand, the rank of $\tf$ restricted to $W'$ equals~$1$
identically (if $\tf=0$ at a point $(a,U)$ then $U\subset\Ker f_a$,
hence $a$ is a singular point of $Y$, see~proposition~\ref{sm}).
Hence, the cokernel of $\tf$ is a line bundle on $W'$, denote 
it by $\CL_{W'}$. The composition of the canonical projection
$V^*\otimes\CO_{W'}\to (\CO_Y\boxtimes\CU^*)_{|W'}$
and of the cokernel morphism $(\CO_Y\boxtimes\CU^*)_{|W'} \to \CL_{W'}$
specifies a map $W'\to Y\times X\times\PP(V)$.
Furthermore, since this morphism factors through
$(\CO_Y\boxtimes\CU^*)_{|W'}$, the map factors through $Y\times\PP_X(\CU)$.
Similarly, it is easy to show that 
the morphism $V^*\otimes\CO_{W'}\to\CL_{W'}$
factors through $(E\boxtimes\CO_X)_{|W'}$ (see the above diagram),
hence the map $W'\to Y\times X\times\PP(V)$
factors through $\PP_Y(E^*)\times X$.
Therefore, we obtain a pair of maps $W'\to \PP_Y(E^*)$ and $W'\to\PP_X(\CU)$,
such that the compositions $W'\to \PP_Y(E^*) \to Q\subset\PP(V)$ and
$W'\to\PP_X(\CU) \to Q\subset\PP(V)$ coincide. Thus we obtain a map
$W'\to\PP_Y(E^*)\times_Q\PP_X(\CU) = W$. It is easy to see that
this map is inverse to the map $j$ above.

Thus we have proved that $j:W\to W'$ is a scheme-theoretical isomorphism.
Moreover, the above arguments show that the cokernel of $\tf$, $\CL_{W'}$,
is isomorphic to $j_*\phipsi^*\CO_Q(1)$. Therefore the sequence~(\ref{jw})
is exact at least at the right two terms. Furthermore, the above arguments
also prove exactness at the left term. It remains to check that the
embedding $E^*(-1) \boxtimes \CO_X \to \Ker\tf$ is an isomorphism.
Indeed, this is true because both sheaves are reflexive of rank~$2$ and
\begin{multline*}
\det\Ker\tf \cong
\det(\CO_Y(-1)\boxtimes V/\CU)\otimes\det(\CO_Y\boxtimes \CU^*)^{-1}\otimes
\det j_*\phipsi^*\CO_Q(1)
\cong \\ \cong
(\CO_Y(-4)\boxtimes\CO_X(1))\otimes(\CO_Y\boxtimes\CO_X(1))^{-1} \cong
\CO_Y(-4)\boxtimes\CO_X \cong
\det(E^*(-1) \boxtimes \CO_X).
\end{multline*}
Here $\det j_*\phipsi^*\CO_Q(1) \cong \CO_{Y\times X}$ 
because $\codim\supp(j_*\phipsi^*\CO_Q(1)) = 2$.

$(ii)$
%Consider $\PP(A)\times\PP_X(\CU)$.
Let $\hf$ denote the composition of the following morphisms 
on $\PP(A)\times\PP_X(\CU)$:
$$
\CO_{\PP(A)}(-1)\boxtimes p_X^*V/\CU \to
\CO_{\PP(A)}\boxtimes p_X^*\CU^* \to
\CO_{\PP(A)}\boxtimes \phi^*\CO_Q(1),
$$
where the first morphism is defined similarly to the morphism $\tf$ 
in $(i)$, and the second morphism is the canonical projection.
Let $W''\subset\PP(A)\times\PP_X(\CU)$ denote the zero scheme of $\hf$.
In the other words $W''$ is the zero scheme of a section of the vector bundle
$\CO_{\PP(A)}(1)\boxtimes(\phi^*\CO_Q(1)\otimes p_X^*(V/\CU)^*)$,
corresponding to $\hf$. We are going to prove that the map
$\PP(A)\times\PP_X(\CU) \exto{\id\times p_X} \PP(A)\times X$
induces an isomorphism of $W''\subset\PP(A)\times\PP_X(\CU)$
to $W'\subset Y\times X\subset\PP(A)\times X$.

Indeed, the definition of $\hf$ shows that the pullback under
$\id\times p_X$ of $(\alpha\times\id_X)_*\tf$ degenerates on $W''$, hence
$(\id\times p_X)(W'')\subset W'$. Similarly, the map
$(\CO_Y\boxtimes\CU^*)_{|W'} \to \phipsi^*\CO_Q(1)$ from~(\ref{jw})
specifies an embedding $W'\to Y\times\PP_X(\CU)\subset\PP(A)\times\PP_X(\CU)$,
and it is clear that the pullback of $\hf$ under this embedding vanishes.
Thus we obtain the iverse map $W'\to W''$.

Further, it is easy to see that the composition of $j:W\to W'$,
and of the above isomorphism $W'\to W''$ coincides with $\lambda$.
Finally,
$\codim W''=
\rank(\CO_{\PP(A)}(1)\boxtimes(\phi^*\CO_Q(1)\otimes p_X^*(V/\CU)^*))$,
hence the structure sheaf $\lambda_*\CO_W$ admits a Koszul
resolution~(\ref{jw1}).
\end{proof}

\section{Derived categories}

\subsection*{Preliminaries}

Let $\D$ be a triangulated category~\cite{V,GM}.
An important example of a triangulated category 
is $\D^b(M)$, the bounded derived category 
of coherent sheaves on a smooth projective 
variety $M$. We briefly remind some definitions 
and results from~\cite{BK,B,BO,Or} and \cite{Br}.

\begin{definition}[\cite{B}]
An object $F\in\D$ is called {\em exceptional}\/ if $\Hom(F,F)=\C$
and $\Ext^p(F,F)=0$ for $p\ne 0$. A collection of exceptional
objects $(F_1,\dots,F_k)$ is called {\em exceptional}\/ if
$\Ext^p(F_i,F_j)=0$ for $i>j$ and all $p\in\ZZ$.
\end{definition}

%Let $\CA$ be a strictly full triangulated subcategory of $\D$
%and $\alpha:\CA\to\D$ be the embedding functor.

\begin{definition}[\cite{B}]
A strictly full triangulated subcategory $\CA\subset\D$ 
is {\em admissible}\/ if the embedding functor $\CA\to\D$ 
admits the left and the right adjoint functors $\D\to\CA$.
\end{definition}

\begin{proposition}[\cite{B}]\label{ex_adm}
Let $(F_1,\dots,F_k)$ be an exceptional collection in $\D$.
The triangulated subcategory $\lan F_1,\dots,F_k\ran \subset \D$ 
generated by objects $F_1$, \dots, $F_k$ is admissible.
\end{proposition}

If $\CA$ is a full triangulated subcategory of $\D$
then the {\em right orthogonal}\/ to $\CA$ in $\D$ is the full
subcategory $\CA^\perp\subset\D$ consisting of all objects
$G\in\D$ such that $\Hom(F,G)=0$ for all $F\in\CA$.
Similarly, the {\em left orthogonal}\/ to $\CA$ in $\D$ is the full
subcategory ${}^\perp\!\CA\subset\D$ consisting of all objects
$G\in\D$ such that $\Hom(G,F)=0$ for all $F\in\CA$.

\begin{definition}[\cite{BO}]
A sequence of admissible subcategories $(\CA_1,\dots,\CA_n)$ in $\D$
is {\em semiorthogonal}\/ if $\CA_j\subset\CA_i^\perp$ for $i>j$.
Triangulated subcategory of $\D$ generated by subcategoires 
$\CA_1$, \dots, $\CA_n$ is denoted by $\lan\CA_1,\dots,\CA_n\ran$.
A semiorthogonal collection $(\CA_1,\dots,\CA_n)$ is {\em full}\/ 
if $\lan\CA_1,\dots,\CA_n\ran = \D$. A~full semiorthogonal collection 
in $\D$ is called a {\em semiorthogonal decomposition}\/ of $\D$.
\end{definition}

\begin{definition}[\cite{BK}]
A covariant additive functor $\FS_\D:\D\to\D$ is a {\em Serre functor}\/
if it is a category equivalence and for all objects $F,G\in\D$ there are given 
bi-functorial isomorphisms $\varphi_{F,G}:\Hom(F,G)\to\Hom(G,\FS_\D(F))^*$
such that the composition
$$
(\varphi^{-1}_{G,\FS_\D(F)})^*\circ\varphi_{F,G}:
\Hom(F,G) \to \Hom(G,\FS_\D(F))^* \to \Hom(\FS_\D(F),\FS_\D(G))
$$
coinsides with the isomorphism induced by $\FS_\D$.
\end{definition}

\begin{proposition}[\cite{BK}]\label{serre_bk}
If a Serre functor exists then it is unique up to a canonical
functorial isomorphism. If $\D = \D^b(M)$ then
$\FS_\D(F):=F\otimes\omega_M[\dim M]$ is a Serre functor.
\end{proposition}

\begin{proposition}[\cite{B}]\label{sos_sod}
If $\D$ admits a Serre functor and $(\CA_1,\dots,\CA_n)$
is a semiorthogonal sequence of admissible subcategories,
then $\D=\lan\CA_0,\CA_1,\dots,\CA_n\ran$ and 
$\D=\lan\CA_1,\dots,\CA_n,\CA_{n+1}\ran$
are semiotrhogonal decompositions, where
$\CA_0=\lan\CA_1,\dots,\CA_n\ran^\perp$ and
$\CA_{n+1}={}^\perp\lan\CA_1,\dots,\CA_n\ran$.
\end{proposition}

\begin{proposition}[\cite{B}]\label{mut_funct}
If $\D$ admits a Serre functor and $\CA\subset\D$ is admissible 
then there exist exact functors $L_{\CA}:\D\to\CA^\perp$ and
$R_{\CA}:\D\to{}^\perp\!\CA$ inducing equivalences
${}^\perp\!\CA\to\CA^\perp$, $\CA^\perp\to{}^\perp\!\CA$,
such that $L_\CA(\CA) = 0$, $R_\CA(\CA) = 0$,
$(L_{\CA})_{|{}^\perp\!\CA} = \FS_{\D}\circ\FS_{{}^\perp\!\CA}^{-1}$, and
$(R_{\CA})_{|\CA^\perp} = \FS_{\D}^{-1}\circ\FS_{\CA^\perp}$.
Moreover, such functors are unique up to a canonical functorial isomorphism.
\end{proposition}

\begin{proposition}[\cite{B,BO}]\label{translations}
Let $\D=\lan\CA_1,\dots,\CA_n\ran$ be a semiorthogonal decmposition. 
If $\D$ admits a Serre functor then for any $1\le k\le n-1$
we have semiorthogonal decompositions
$$
\begin{array}{rcl}
\D &=& \lan\CA_1,\dots,\CA_{k-1},\CA_{k+1},
R_{\CA_{k+1}}\CA_k,\CA_{k+2},\dots,\CA_n\ran,\\
\D &=& \lan\CA_1,\dots,\CA_{k-1},L_{\CA_{k}}\CA_{k+1},
\CA_k,\CA_{k+2},\dots,\CA_n\ran,
\end{array}
$$
and $R_{\CA_{k+1}}:\CA_k \to R_{\CA_{k+1}}\CA_k$,
$L_{\CA_k}:\CA_{k+1} \to L_{\CA_k}\CA_{k+1}$ are equivalences.
If additionally $\CA_{k+1}\subset\CA_k^\perp$ {\rm(}i.e.\ $\CA_k$
and $\CA_{k+1}$ are completely orthogonal{\rm)}, then
$L_{\CA_k}\CA_{k+1} = \CA_{k+1}$, $R_{\CA_{k+1}}\CA_k = \CA_k$.
\end{proposition}

We will call these operations on semiorthogonal decompositions
the (right) {\em mutation of $\CA_k$ through $\CA_{k+1}$}\/ and
the (left) {\em mutation of $\CA_{k+1}$ through $\CA_k$}\/ respectively.
If $\CA=\lan F\ran$ we will denote mutation functors, $L_\CA$ and $R_\CA$, 
by $L_F$ and $R_F$ respectively.

\begin{lemma}\label{philr}
If $\Phi$ is an autoequivalence of $\D$ then 
we have canonical isomorphisms of functors
$\Phi\circ L_\CA \cong L_{\Phi(\CA)}\circ \Phi$,
$\Phi\circ R_\CA \cong R_{\Phi(\CA)}\circ \Phi$.
\end{lemma}

\begin{proposition}[\cite{B}]\label{mutations}
If $\CA_k$ and $\CA_{k+1}$ are generated by exceptional objects
$F_k$ and $F_{k+1}$ respectively, then $L_{\CA_k}\CA_{k+1}$ and
$R_{\CA_{k+1}}\CA_k$ are generated by exceptional objects
$L_{F_k}F_{k+1}$ and $R_{F_{k+1}}F_k$ respectively,
defined by the following exact triangles
$$
\RHom(F_k,F_{k+1})\otimes F_k 
\exto{\ev} F_{k+1} \to
L_{F_k}F_{k+1},
\qquad
R_{F_{k+1}}F_k \to 
F_k \exto{\ev^*} 
\RHom(F_k,F_{k+1})^*\otimes F_{k+1},
$$
where $\ev$ and $\ev^*$ denote the canonical evaluation and coevaluation
homomorphisms.
\end{proposition}

Let $M_1$, $M_2$ be smooth projective varieties and
let $p_i:M_1\times M_2 \to M_i$ denote the projections.
Take any $K\in\D^b(M_1\times M_2)$ and define
$\Phi_K(F) := {p_2}_*(p_1^*F\otimes K)$.
Then $\Phi_K$ is an exact functor $\D^b(M_1)\to\D^b(M_2)$,
the {\em kernel functor}\/ with kernel $K$. Kernel functors
can be thought of as analogues of correspondences
on categorical level.

\begin{lemma}\label{comp_kern2}
If $K\in\D^b(M_2\times M_3)$, $F_1\in\D^b(M_1)$, $F_2\in\D^b(M_2)$, then 
$\Phi_K\cdot\Phi_{F_1\boxtimes F_2} \cong \Phi_{F_1\boxtimes\Phi_K(F_2)}$.
\end{lemma}

\begin{proposition}\label{mut_kernels}
If $M$ is a smooth projective variety, $\D=\D^b(M)$,
and $F\in\D$ is an exceptional object then the mutation 
functors $L_F$, $R_F$ are kernel functors
given by the kernels ${}_F\SK$ and $\SK_F$ on $M\times M$ defined
by the following exact triangles
$$
\RCHom(F,\CO_M)\boxtimes F \exto{\ev} \Delta_*\CO_M \to {}_F\SK,\qquad
\SK_F \to \Delta_*\CO_M \exto{\ev^*} \RCHom(F,\omega_M[\dim M])\boxtimes F,
$$
where $\ev$ and $\ev^*$ are the evaluation and coevaluation
homomorphisms, and $\Delta:M\to M\times M$ is the diagonal.
\end{proposition}

Let $M$ be a smooth projective variety and
let $\CE$ be a rank $r$ vector bundle on $M$.
Consider its projectivization $\PP_M(\CE)$
and denote by $p:\PP_M(\CE)\to M$ the projection
and by $L=\CO_{\PP_M(\CE)/M}(1)$ a Grothendieck
relatively ample line bundle.

\begin{proposition}[\cite{Or}]\label{sod_proj}
If $\D^b(M) = \lan\CA_1,\dots,\CA_n\ran$ is a semiorthogonal decomposition then
$$
\D^b(\PP_M(\CE)) =
\lan L^k\otimes p^*\CA_1,\dots,L^k\otimes p^*\CA_n,\dots,
L^{k+r-1}\otimes p^*\CA_1,\dots,L^{k+r-1}\otimes p^*\CA_n\ran
$$
is a semiorthogonal decomposition for any $k\in\ZZ$.
\end{proposition}

We also will need the following reformulation of results
of Bridgeland.

\begin{theorem}[\cite{Br}]\label{flop_equ}
Let $M$ be a smooth projective variety and let $\psi:M\to M'$ be
a crepant contraction of relative dimension~$1$. 
Let $Z\subset M$ denote the exceptional locus of $\psi$.
Assume that $\psi^+:M^+\to M'$ is a flop of $\psi$ with $M^+$ smooth
and let $Z^+\subset M^+$ denote the exceptional locus of $\psi^+$,
so that $\psi^{-1}\cdot\psi^+:M^+-Z^+ \to M-Z$ is an isomorphism.
For any point $s\in M^+$ let $j_s:M\to M\times M^+$ denote
the corresponding embedding.
If $K\in\D^b(M\times M^+)$ is an object, such that
for any point $s\in M^+$ we have either
\begin{itemize}
\item $Lj_s^*K \cong \CO_{\psi^{-1}\cdot\psi^+(s)}$, if $s\in M^+-Z^+$;
\item we have an exact triangle
$Lj_s^*K \to \CO_L \exto{\epsilon} \CO_L(-1)[2]$ with $\epsilon\ne0$,
where $L = \psi^{-1}\cdot\psi^+(s) \cong \PP^1$, if $s\in Z^+$.
\end{itemize}
Then the kernel functor $\Phi_K:\D^b(M)\to\D^b(M^+)$
is an equivalence.
\end{theorem}

\subsection*{Derived categories of $Y$ and $X$}

Let $Y\subset\PP(A)$ be a smooth cubic threefold and
let $X$ be a smooth $V_{14}$ threefold.
To avoid an abuse of notation, let us denote
by $\CO(y)$ the sheaf $\CO_{\PP(A)}(1)$
and its pullbacks to $Y$, $\PP_Y(E^*)$, $W$ etc.,
by $\CO(x)$ the sheaf $\CO_{\Gr(2,V)}(1)$
and its pullbacks to $X$, $\PP_X(\CU)$, $W$ etc.,
and by $\CO(e)$ the sheaf $\CO_{\PP(V)}(1)$
and its pullbacks to $Q$, $\PP_Y(E^*)$, $\PP_X(\CU)$, $W$ etc.

\begin{lemma}\label{ex}
The pairs $(\CO,\CO(y))$ in $\D^b(Y)$ and $(\CO,\CU^*)$ in $\D^b(X)$
are exceptional.
\end{lemma}
\begin{proof}
Straightforward computations using the Koszul resolutions of $Y$ in $\PP(A)$
and of $X$ in $\Gr(2,V)$ and Borel--Bott--Weil theorem.
\end{proof}

The subcategories $\lan\CO,\CO(y)\ran\subset\D^b(Y)$ and
$\lan\CO,\CU^*\ran\subset\D^b(X)$ are admissible 
by proposition~\ref{ex_adm}, hence by proposition~\ref{sos_sod} 
we obtain semiorthogonal decompositions
\begin{equation}\label{sodec}
\D^b(X) = \lan\CO,\CU^*,\CA_X\ran,\qquad
\D^b(Y) = \lan\CA_Y,\CO,\CO(y)\ran,
\end{equation}
where $\CA_X = {}^\perp\lan\CO,\CU^*\ran \subset D^b(X)$
and $\CA_Y = \lan\CO,\CO(y)\ran^\perp \subset \D^b(Y)$.

\begin{theorem}\label{equiv}
If $Y$ is the Pfaffian cubic of $X$ then
categories $\CA_X$ and $\CA_Y$ are equivalent.
\end{theorem}

\begin{corollary}
If $X$ and $X'$ are birational then $\CA_X$ and $\CA_{X'}$ are equivalent.
\end{corollary}
\begin{proof}
If $X$ and $X'$ are birational then their Pfaffian cubics $Y$ and $Y'$
are isomorphic by remark~\ref{xy_bir}, hence 
$\CA_X \cong \CA_Y \cong \CA_{Y'} \CA_{X'}$.
\end{proof}

Note that a triangulated category generated by an exceptional object
is equivalent to the derived category of $\C$-vector spaces. 
Therefore we have

\begin{corollary}
If $Y$ is the Pfaffian cubic of $X$ then
derived categories $\D^b(X)$ and $\D^b(Y)$ admit semiorthogonal
decompositions with pairwise equivalent summands.
\end{corollary}

The rest of the section is devoted to the proof of Theorem~\ref{equiv}.
We begin with a short plan of the proof. From now on we assume that
$Y$ is the Pfaffian cubic of $X$, $E$ is the corresponding theta-bundle,
so that theorem~\ref{main} holds.

\begin{description}
\item[Step 1]
First of all, we replace for convenience the decomposition~(\ref{sodec}) 
of $\D^b(Y)$ by the decomposition $\D^b(Y) = \lan\CO(-y),\CA_Y,\CO\ran$. 
This is done by mutating $\CO(y)$ to the left, since $\omega_Y\cong\CO(-2y)$.
Further, we note that $\CO(e)$ is a Grothendieck relatively ample line 
bundle both for $\PP_Y(E^*)\to Y$ and for $\PP_X(\CU)\to X$. Hence
by proposition~\ref{sod_proj} we obtain the following
semiorthogonal decompositions
\begin{eqnarray}
\D^b(\PP_X(\CU)) &=& \lan\CO(-e),\CU^*(-e),\CA_X(-e),\CO,\CU^*,\CA_X\ran,
\label{sodectx}\\
\D^b(\PP_Y(E^*)) &=& \lan\CO(-y),\CA_Y,\CO,\CO(e-y),\CA_Y(e),\CO(e)\ran,
\label{sodecty}
\end{eqnarray}
where $p_X^*$ and $p_Y^*$ are omitted for brevity.

\item[Step 2]
We perform with the decomposition $(\ref{sodecty})$ a sequence
of mutations (described below) and obtain the following semiorthogonal 
decomposition
\begin{equation}\label{sodectynew}
\begin{array}{rlll}
\D^b(\PP_Y(E^*)) = \lan
& \CO(-e), & L_\CO\CO(2e-y),      & R_{\CO(2e-y)}\CA_Y(e), \\
& \CO,     & L_{\CO(e)}\CO(3e-y), & R_{\CO(3e-y)}\CA_Y(2e) \quad \ran.
\end{array}
\end{equation}

\item[Step 3]
Let $K= i_*\CO_W$, where 
$i:W=\PP_Y(E^*) \times_{Q} \PP_X(\CU) \to \PP_Y(E^*) \times \PP_X(\CU)$
is the embedding. We show that the kernel $K$
satisfies the conditions of theorem~\ref{flop_equ}. It follows that
the kernel functor $\Phi_K:\D^b(\PP_Y(E^*)) \to \D^b(\PP_X(\CU))$ 
is an equivalence. We check also that $\Phi_K$ commutes 
with tensoring by pullbacks of sheaves from~$Q$.

\item[Step 4]
We show that
$$
\Phi_K(\CO) \cong \CO,\quad
\Phi_K(L_{\CO}\CO(2e-y)) \cong \CU^*(-e),\quad\text{and}\quad
{p_X}_*\Phi_K(R_{\CO(2e-y)}\CA_Y(e)) = 0.
$$
Lemma~\ref{philr} implies that 
$L_{\CO(e)}\CO(3e-y) \cong (L_\CO\CO(2e-y))\otimes\CO(e)$.
Since $\Phi_K$ commutes with tensoring by $\CO(-e)$ it follows that
$$
\Phi_K(\CO(-e)) \cong \CO(-e),\quad\text{and}\quad
\Phi_K(L_{\CO(e)}\CO(3e-y)) \cong \CU^*.
$$
These isomorphisms show that $\Phi_K$ takes the first line 
of the collection~(\ref{sodectynew}) to the subcategory
$p_X^*\D^b(X)\otimes\CO(-e) =
\lan\CO(-e),\CU^*(-e),\CA_X(-e)\ran \subset \D^b(\PP_X(\CU))$.
Lemma~\ref{philr} implies that the second line of~(\ref{sodectynew}) 
is equal to the first line tensored by $\CO(e)$, therefore $\Phi_K$ 
must induce an equivalence 
$$
\lan \CO, L_{\CO(e)}\CO(3e-y), R_{\CO(3e-y)}\CA_Y(2e) \ran \exto{\Phi_K}
p_X^*\D^b(X) = \lan\CO,\CU^*,\CA_X\ran \subset \D^b(\PP_X(\CU)).
$$
Finally, since $\Phi_K(\CO)\cong\CO$ and 
$\Phi_K(L_{\CO(e)}\CO(3e-y)) \cong \CU^*$ 
it follows that $\Phi_K$ must induce an equivalence 
$R_{\CO(3e-y)}\CA_Y(2e) \to \CA_X \subset \D^b(\PP_X(\CU))$.
Summarizing, we see that
\begin{equation}\label{exp_equ}
\Phi(A) = {p_X}_*(\Phi_K(R_{\CO(3e-y)}(p_Y^*(A)\otimes\CO(2e))))
\end{equation}
is an equivalence $\CA_Y \to \CA_X$.
\end{description}

Now we start implementing above steps. Step~1 is already
quite clear, so we can pass to Step~2.

\subsection*{Step 2}

First of all, we note that $\omega_{\PP_Y(E^*)} \cong \CO(-2e)$
(see the proof of theorem~\ref{flop}). Further, we will need
the following

\begin{lemma}\label{ext_ty}
In $\D^b(\PP_Y(E^*))$ we have

$(i)$ $\Ext^p(\CO,\CO(e-y))=0$ for all $p\in\ZZ$.

$(ii)$ $\Ext^p(\CO,\CO(2e-y))=
\begin{cases}\C, & \text{if $p=1$}\\0, & \text{if $p\ne 1$}\end{cases}$

$(iii)$ $\Ext^p(\CO(-e),R_{\CO(e-y)}F)=0$ 
for any $p\in\ZZ$ and any $F\in\CA_Y$.
\end{lemma}
\begin{proof}
$(i)$
$\Ext^\cdot(\CO,\CO(e-y)) =
H^\cdot(\PP_Y(E^*),\CO(e-y)) =
\HH^\cdot(Y,{p_Y}_*(\CO(e-y))) =
H^\cdot(Y,E(-1)) = 0$ 
by theorem~\ref{v14}.

$(ii)$ Similarly, we have
$\Ext^\cdot(\CO,\CO(2e-y)) =
\HH^\cdot(Y,{p_Y}_*(\CO(2e-y))) =
H^\cdot(Y,S^2E(-1))$
and it remains to apply lemma~\ref{temy}.

$(iii)$ Using lemma~\ref{philr} and theorem~\ref{phi_k_equ} we deduce that
$$
\Ext^\cdot(\CO(-e),R_{\CO(e-y)}F) \cong
\Ext^\cdot(\CO,R_{\CO(2e-y)}F(e)) \cong 
\Ext^\cdot(\Phi_K(\CO),\Phi_K(R_{\CO(2e-y)}F(e))).
$$
But $\Phi_K(\CO)\cong\CO$ by proposition~\ref{phik1}, hence
$$
\Ext^\cdot(\CO(-e),R_{\CO(e-y)}F) \cong
\HH^\cdot(\PP_X(\CU),\Phi_K(R_{\CO(2e-y)}F(e))) \cong
\HH^\cdot(X,{p_X}_*\Phi_K(R_{\CO(2e-y)}F(e))) 
$$
and it remains to note that
${p_X}_*\Phi_K(R_{\CO(2e-y)}F(e)) = 0$
by proposition~\ref{px_phik_r_py}.
%
%\begin{multline*}
%\Ext^\cdot(\CO(-e),R_{\CO(e-y)}F) \cong
%\Ext^\cdot(\CO,R_{\CO(2e-y)}F(e)) \cong \\ \cong
%\Ext^\cdot(\Phi_K(\CO),\Phi_K(R_{\CO(2e-y)}F(e)) \cong 
%\Ext^\cdot(\CO,\Phi_K(R_{\CO(2e-y)}F(e)) \cong \\ \cong
%H^\cdot(\PP_X(\CU),\Phi_K(R_{\CO(2e-y)}F(e)) \cong
%H^\cdot(X,{p_X}_*\Phi_K(R_{\CO(2e-y)}F(e)) = 0
%\end{multline*}
%by lemma~\ref{philr} and propositions~\ref{phik1} 
%and~\ref{px_phik_r_py}, see below.
\end{proof}

Now, we explain the sequence of transformations.
We start with semiortogonal deocmposition
$\D^b(\PP_Y(E^*)) = \lan\CO(-y),\CA_Y,\CO,\CO(e-y),\CA_Y(e),\CO(e)\ran$.

\begin{enumerate}
\item
We mutate $\CO(-y)$ to the right; it is get twisted by
$\CO(2e)$, the anticanonical class of $\PP_Y(E^*)$:
$$
\D^b(\PP_Y(E^*)) = \lan\CA_Y,\CO,\CO(e-y),\CA_Y(e),\CO(e),\CO(2e-y)\ran.
$$
\item
We mutate $\CO$ through $\CO(e-y)$ and
$\CO(e)$ through $\CO(2e-y)$; lemma~\ref{ext_ty}~(i) and
proposition~\ref{translations} imply that
$R_{\CO(e-y)}\CO = \CO$, $R_{\CO(2e-y)}\CO(e) = \CO(e)$,
and we get
$$
\D^b(\PP_Y(E^*)) = \lan\CA_Y,\CO(e-y),\CO,\CA_Y(e),\CO(2e-y),\CO(e)\ran.
$$
\item
We mutate $\CA_Y$ through $\CO(e-y)$ and $\CA_Y(e)$ through $\CO(2e-y)$:
$$
\D^b(\PP_Y(E^*)) =
\lan\CO(e-y),R_{\CO(e-y)}\CA_Y,\CO,\CO(2e-y),R_{\CO(2e-y)}\CA_Y(e),\CO(e)\ran.
$$
\item
We mutate $\CO(e-y)$ to the right; it is get twisted by $\CO(2e)$:
$$
\D^b(\PP_Y(E^*)) =
\lan R_{\CO(e-y)}\CA_Y,\CO,\CO(2e-y),R_{\CO(2e-y)}\CA_Y(e),\CO(e),\CO(3e-y)\ran.
$$
\item
We mutate $\CO(2e-y)$ through $\CO$ and $\CO(3e-y)$ through $\CO(e)$;
lemma~\ref{ext_ty}~(ii) and proposition~\ref{mutations} imply that
$L_{\CO}\CO(2e-y)$ and $L_{\CO(e)}\CO(3e-y)$ are the unique 
nontrivial extensions
\begin{equation}\label{lo}
\begin{array}{lllllllll}
0 &\to & \CO(2e-y) &\to & L_{\CO}\CO(2e-y) &\to & \CO &\to & 0,\\
0 &\to & \CO(3e-y) &\to & L_{\CO(e)}\CO(3e-y) &\to & \CO(e) &\to & 0,
\end{array}
\end{equation}
and we get:
$$
\D^b(\PP_Y(E^*)) =
\lan R_{\CO(e-y)}\CA_Y,L_{\CO}\CO(2e-y),\CO,
R_{\CO(2e-y)}\CA_Y(e),L_{\CO(e)}\CO(3e-y),\CO(e)\ran.
$$
\item
We mutate $\CO(e)$ to the left; it is get twisted by $\CO(-2e)$:
$$
\D^b(\PP_Y(E^*)) =
\lan\CO(-e),R_{\CO(e-y)}\CA_Y,L_{\CO}\CO(2e-y),
\CO,R_{\CO(2e-y)}\CA_Y(e),L_{\CO(e)}\CO(3e-y)\ran.
$$
\item
We mutate $\CO(-e)$ through $R_{\CO(e-y)}\CA_Y$ and
$\CO$ through $R_{\CO(2e-y)}\CA_Y(e)$; lemma~\ref{ext_ty}~(iii) 
and proposition~\ref{translations} imply that
the mutations coincide with transpositions:
$$
\D^b(\PP_Y(E^*)) =
\lan R_{\CO(e-y)}\CA_Y,\CO(-e),L_{\CO}\CO(2e-y),
R_{\CO(2e-y)}\CA_Y(e),\CO,L_{\CO(e)}\CO(3e-y)\ran.
$$
\item
We mutate $R_{\CO(e-y)}\CA_Y$ to the right;
it is get twisted by $\CO(2e)$ and once again
using lemma~\ref{philr} we get the desired 
decomposition~(\ref{sodectynew}).
\end{enumerate}

This completes Step~2.

\subsection*{Step 3}

We adopt the notation of proposition~\ref{fp} and of theorem~\ref{flop_equ}.

\begin{theorem}\label{phi_k_equ}
If $K = i_*\CO_W$ then the kernel functor 
$\Phi_K:\D^b(\PP_Y(E^*)) \to \D^b(\PP_X(\CU))$
is an equivalence. Moreover, $\Phi_K$ commutes with tensoring
by pullbacks of bundles from $Q$.
\end{theorem}
\begin{proof}
We must check the conditions of theorem~\ref{flop_equ}.
Take an arbitrary point $s\in\PP_X(\CU)$. Then
\begin{multline*}
(\alpha p_Y\times\id)_*{j_s}_*j_s^* K \cong
(\alpha p_Y\times\id)_*{j_s}_*j_s^* i_*\CO_W \cong
(\alpha p_Y\times\id)_*\big(i_*\CO_W \otimes {j_s}_*\CO_{\PP_Y(E^*)}\big) 
\cong \\ \cong
(\alpha p_Y\times\id)_*i_*i^* {j_s}_*\CO_{\PP_Y(E^*)} \cong
\lambda_*i^* \pi_2^*\CO_s \cong
\lambda_*\lambda^*{\pi'_2}^*\CO_s \cong
{\pi'_2}^*\CO_s \otimes \lambda_*\CO_W.
\end{multline*}
Here $\pi_2$ and $\pi_2'$ denote the projections
of $\PP_Y(E^*)\times\PP_X(\CU)$ and
$\PP(A)\times\PP_X(\CU)$ to $\PP_X(\CU)$
and $p_Y$ is considered as a map $\PP_Y(E^*)\to\PP(A)$:
$$
\xymatrix@R=10pt{
&&& W \arrow[dl]_-{i} \arrow[dr]^-{\lambda} \\
\PP_Y(E^*) \arrow[rr]^-{j_s} &&
\PP_Y(E^*)\times\PP_X(\CU) 
\arrow[rr]^{\alpha p_Y\times\id} \arrow[dr]_{\pi_2} &&
\PP(A)\times\PP_X(\CU) \arrow[dl]^{\pi_2'} \\
&&& \PP_X(\CU)
}
$$
The RHS of the above chain of isomorphisms can be computed with a help 
of resolution~(\ref{jw1}). It follows, that for any point 
$s\not\in\SX\subset\PP_X(\CU)$ it is isomorphic to the 
structure sheaf of the point $p_Y(\theta^{-1}(s))\in\PP(A)$,
while for $s\in\SX\subset\PP_X(\CU)$ we get an exact triangle
$$
{\pi'_2}^*\CO_s \otimes \lambda_*\CO_W
\to \CO_{M\times s} \to \CO_{M\times s}(-1)[2],
$$
where $M=p_Y(\psi^{-1}(\phi(s)))\subset\PP(A)$.
Since the object $j_s^*K$ is supported on $\psi^{-1}(\phi(s))$ and
by proposition~\ref{fp}~$(ii)$ the map $p_Y:\psi^{-1}(\phi(s)) \to \PP(A)$ 
is a closed embedding, it follows that we have an exact triangle
$$
j_s^*K \to \CO_\TM \exto{\epsilon} \CO_\TM(-1)[2],
$$
where $\TM=\psi^{-1}(\phi(s)) \subset \PP_Y(E^*)$.
It remains to check that $\epsilon\ne 0$. Note that $\epsilon=0$
would imply $j_s^*K\cong \CO_\TM \oplus \CO_\TM(-1)[1]$, hence
$\Hom^1(j_s^*K,\CO_\TM(-1))\ne0$. Thus, it suffices
to check that $\Hom^1(j_s^*K,\CO_\TM(-1)) = 0$.
To this end we consider the following diagram
$$
\xymatrix{
\PP_Y(E^*) \arrow[r]^-{j_s} \arrow[d]_{\psi} &
\PP_Y(E^*)\times\PP_X(\CU) \arrow[d]_{\psi\times\phi} &
W \arrow[l]_-{i} \arrow[d]_{\phipsi} \\
Q \arrow[r]^-{j_c} &
Q\times Q &
Q \arrow[l]_-{\Delta}
}
$$
where $c=\phi(s)\in Q$, $j_c(v)=(v,c)\in Q\times Q$
and $\Delta$ is the diagonal.
In this diagram the right square is Cartesian
and $\Delta$ is a closed embedding, hence lemma~\ref{my}
implies that there is a functorial morphism
$(\psi\times\phi)^*\Delta_* \to i_*\phipsi^*$, and furthermore,
for any $F\in\D^{\le0}(Q)$ the object $F'$ in the exact triangle
$$
F' \to (\psi\times\phi)^*\Delta_* F \to i_*\phipsi^* F
$$
is contained in $\D^{\le -1}(\PP_Y(E^*)\times\PP_X(\CU))$.
Applying $j_s^*$ and using $j_s^*(\psi\times\phi)^*=\psi^*j_c^*$
we get an exact triangle
$$
F'' \to \psi^*j_c^*\Delta_* F \to j_s^*i_*\phipsi^* F,
$$
with $F'' = j_s^*F'\in\D^{\le-1}(\PP_Y(E^*))$ since $j_s^*$ is right exact.
Substituting $F=\CO_Q$ and using isomorphisms $\phipsi^*\CO_Q\cong\CO_W$,
%$j_c^*\Delta_*\CO_Q \cong \CO_c$, 
we obtain a triangle
$$
F'' \to \psi^*j_c^*\Delta_*\CO_Q \to j_s^*K.
$$
Applying the functor $\Hom(-,\CO_\TM(-1))$ and using
$$
\Hom^\cdot(\psi^*j_c^*\Delta_*\CO_Q,\CO_\TM(-1)) \cong
\Hom^\cdot(j_c^*\Delta_*\CO_Q,\psi_*\CO_\TM(-1)) =
\Hom^\cdot(j_c^*\Delta_*\CO_Q,0) = 0,
$$
we deduce that $\Hom^1(j_s^*K,\CO_\TM(-1)) = \Hom(F'',\CO_\TM(-1)) = 0$,
since we have $F''\in\D^{\le-1}(\PP_Y(E^*))$ and
$\CO_\TM(-1)\in\D^{\ge0}(\PP_Y(E^*))$.

Now theorem~\ref{flop_equ} implies that the functor 
$\Phi_K:\D^b(\PP_Y(E^*)) \to \D^b((\PP_X(\CU))$
is an equivalence.

Finally, let $\CV$ be an arbitrary vector bundle on $Q$. Then 
the functor $F\mapsto \Phi_K(F\otimes\psi^*\CV)$ is a kernel functor
with kernel 
$K\otimes\pi_1^*\psi^*\CV =
i_*\CO_W\otimes\pi_1^*\psi^*\CV =
i_*i^*\pi_1^*\psi^*\CV =
i_*\phipsi^*\CV$,
and the functor $F\mapsto \Phi_K(F)\otimes\phi^*\CV$ is a kernel functor
with kernel
$K\otimes\pi_2^*\phi^*\CV =
i_*\CO_W\otimes\pi_2^*\phi^*\CV =
i_*i^*\pi_2^*\phi^*\CV =
i_*\phipsi^*\CV$,
where $\pi_1$ and $\pi_2$ are projections of $\PP_Y(E^*)\times\PP_X(\CU)$
to the factors. The kernels are isomorphic, hence 
the functors are isomorphic as well.
\end{proof}

\begin{lemma}[cf.\ \cite{Sw}]\label{my}
For any Carthesian square
$$
\xymatrix{
T \arrow[d]_g & 
T' \arrow[l]_{f'} \arrow[d]^{g'} \\ 
S & 
S' \arrow[l]_{f}
}
$$
there is a canonical morphism of functors
$g^*f_* \to f'_*{g'}^*$. Further, if $f$ is affine
then for any $F\in\D^{\le 0}(S')$ we have $F'\in\D^{\le-1}(T)$,
where $F'$ fits into the triangle $F' \to g^*f_*F \to f'_*{g'}^*F$.
\end{lemma}
\begin{proof}
Using the adjunction morphisms for $g$ and $g'$, and an isomorphis
$f_*g'_*\cong g_*f'_*$ we define the morphism of functors as 
the following composition
$$
g^*f_* \to g^*f_*g'_*{g'}^* \to g^*g_*f'_*{g'}^* \to f'_*{g'}^*.
$$
Now, assume that $f$ is affine and let us check that
$F' \in \D^{\le-1}(T)$. The property is local, so we can assume
that $S$ and $T$ are affine, say $S = \Spec A$, $T = \Spec B$.
Then $S' = \Spec A'$, $T' = \Spec B\otimes_A A'$, where $A'$
is a finitely generated $A$-algebra. Note that
$$
f'_*{g'}^* (A') \cong B\otimes_A A',\qquad\text{and}\qquad
g^*f_* (A') \in \D^{\le0}(T),\quad
H^0(g^*f_* (')) \cong B\otimes_A A'.
$$
Taking a resolution of $F\in\D^{\le0}(S')$ by free $A'$-modules
we deduce the claim.
\end{proof}

\subsection*{Step 4}

\begin{proposition}\label{phik1}
We have

$(i)$ $\Phi_K(\CO) = \CO$, $\Phi_K(\CO(-e)) = \CO(-e)$;

$(ii)$ $\Phi_K(L_\CO \CO(2e-y)) = \CU^*(-e)$,
$\Phi_K(L_{\CO(e)} \CO(3e-y)) = \CU^*$.
\end{proposition}
\begin{proof}
First of all we note that for any $F\in\D^b(\PP(A))$ we have
\begin{multline*}
\Phi_K(p_Y^*\alpha^*F) =
{\pi_2}_*(\pi_1^*p_Y^*\alpha^*F \otimes i_*\CO_W) \cong
{\pi_2}_*i_*i^*\pi_1^*p_Y^*\alpha^*F \cong \\ \cong
\xi_*\eta^*p_Y^*\alpha^*F =
{\pi'_2}_*\lambda_*\lambda^*{\pi'_1}^*F =
{\pi'_2}_*({\pi'_1}^*F\otimes\lambda_*\CO_W),
\end{multline*}
where $\pi_1$, $\pi_2$ are the projections of $\PP_Y(E^*)\times\PP_X(\CU)$
to the factors, and 
$\pi'_1$, $\pi'_2$ are the projections of $\PP(A)\times\PP_X(\CU)$
to the factors.

$(i)$ Taking $F=\CO_{\PP(A)}$ and applying~(\ref{jw1}) we get
$\Phi_K(\CO) = \CO$. Further, since $\CO(-e)$ is a pullback of
a line bundle from $Q$, it follows from theorem~\ref{phi_k_equ}
that $\Phi_K(\CO(-e))=\CO(-e)$.

$(ii)$ Taking $F=\CO_{\PP(A)}(-1)$ and applying~(\ref{jw1}) we get
$\Phi_K(\CO(-y)) = R^4{\pi'_2}_*\Lambda^4(V/U)(-4e-5y) = \CO(x-4e)$.
Further, since $\CO(2e)$ is a pullback of a line bundle from $Q$,
it follows from theorem~\ref{phi_k_equ} that 
\begin{equation}\label{xmte}
\Phi_K(\CO(2e-y)) \cong \CO(x-2e).
\end{equation}
Since $L_\CO\CO(2e-y)$ is the unique nontrivial extension of $\CO$
by $\CO(2e-y)$ and since $\Phi_K$ is an equivalence, it follows
that $\Phi_K(L_\CO\CO(2e-y))$ is the unique nontrivial extension
of $\CO$ by $\CO(x-2e)$. On the other hand, it is clear that $\CU^*(-e)$
is such an extension. Hence $\Phi_K(L_\CO\CO(2e-y))\cong\CU^*(-e)$.
Finally, by lemma~\ref{philr} we have 
$L_{\CO(e)}\CO(3e-y) \cong (L_\CO\CO(2e-y)) \otimes \CO(e)$
hence by theorem~\ref{phi_k_equ} $\Phi_K(L_{\CO(e)}\CO(3e-y))\cong\CU^*$. 
\end{proof}

\begin{proposition}\label{px_phik_r_py}
We have ${p_X}_*\Phi_K(R_{\CO(2e-y)}\CA_Y(e)) = 0$.
\end{proposition}
\begin{proof}
First of all, lemma~\ref{philr} implies that
$$
\Phi_K(R_{\CO(2e-y)}\CA_Y(e)) =
\Phi_{K(e)}(R_{\CO(e-y)}\CA_Y) =
(\Phi_{K(e)}\cdot \Phi_F)(\CA_Y),
$$
where $F = \SK_{\CO(e-y)}\in\D^b(\PP_Y(E^*)\times\PP_Y(E^*))$ 
is defined from the following exact triangle 
(cf.~proposition~\ref{mut_kernels})
\begin{multline*}
\{F \to 
\Delta_*\CO_{\PP_Y(E^*)} \exto{\rho}
\RCHom(\CO(e-y),\CO(-2e)[4])\boxtimes\CO(e-y)\} 
= \\ =
\{F \to
\Delta_*\CO_{\PP_Y(E^*)} \exto{\rho} 
\CO(y-3e)\boxtimes\CO(e-y)[4]\}
\end{multline*}
with $\rho\ne0$.
Since the kernel $\Delta_*\CO_{\PP_Y(E^*)}$ gives the identity functor,
it follows from lemma~\ref{comp_kern2} that the functor 
$\Phi_{K(e)}\cdot \Phi_F$ is given by
the kernel $K'\in\D^b(\PP_Y(E^*)\times\PP_X(\CU))$, defined
from the exact triangle
$$
K' \to K(e) \exto{\rho'} \CO(y-3e)\boxtimes\Phi_{K(e)}(\CO(e-y))[4],
$$
where $\rho'=\Phi_{K(e)}(\rho)\ne 0$, since $\Phi_{K(e)}$ is an equivalence.
Further, applying~(\ref{xmte}) we obtain
$\Phi_{K(e)}(\CO(e-y)) = \Phi_K(\CO(2e-y)) = \CO(x-2e)$,
and it follows that we have the following exact triangle:
$$
K' \to K(e) \exto{\rho'} \CO(y-3e)\boxtimes\CO(x-2e)[4].
$$
It is clear that we have
$\Phi_{K'} \cdot p_Y^* = \Phi_{(p_Y\times \id)_*K'}$.
Applying $(p_Y\times \id)_*$ to the above triangle we see that
the resulting functor is $\Phi_{K''}$ where $K''\in\D^b(Y\times \PP_X(\CU))$
is defined from the following exact triangle:
$$
K'' \to (p_Y\times \id)_*K(e) \exto{\rho''}
{p_Y}_*\CO(y-3e)\boxtimes\CO(x-2e)[4].
$$
Since ${p_Y}_*\CO(y-3e) \cong E^*(-y)[-1]$, we have
$$
K'' \to (p_Y\times \id)_*K(e) \exto{\rho''}
E^*(-y)\boxtimes\CO(x-2e)[3].
$$
Let us check that $\rho''\ne0$. Indeed, if $\rho''=0$ then
$K''$ would have $E^*(-y)\boxtimes\CO(x-2e)[2]$ as a direct summand,
hence for any $F\in\CA_Y$, such that $H^\cdot(Y,E^*(-y)\otimes F)\ne 0$
(e.g.\ $F=E(-y)$) the object $\Phi_{K''}(F)$ would have a shift of $\CO(x-2e)$ 
as a direct summand, hence using~(\ref{xmte}) we would obtain
\begin{multline*}
0 \ne
\Hom^\cdot(\Phi_{K''}(F),\CO(x-2e)) =
\Hom^\cdot(\Phi_K(R_{\CO(2e-y)}F(e)),\CO(x-2e)) = \\ =
\Hom^\cdot(\Phi_K(R_{\CO(2e-y)}F(e)),\Phi_K(\CO(2e-y))) =
\Hom^\cdot(R_{\CO(2e-y)}F(e),\CO(2e-y)),
\end{multline*}
which would give a contradiction with proposition~\ref{mut_funct}. 

Thus $\rho''\ne0$. Further, it is clear that we have
${p_X}_*\Phi_{K''} = \Phi_{(\id\times p_X)_*K''}$,
and applying $(\id\times p_X)_*$ to the above triangle we see that
the resulting functor is $\Phi_{K'''}$ where
$K'''\in\D^b(Y\times X)$ is defined from the following exact triangle:
$$
K''' \to (p_Y\times p_X)_*K(e) \exto{\rho'''}
E^*(-y)\boxtimes{p_X}_*(\CO(x-2e))[3].
$$
Since ${p_X}_*\CO(x-2e) \cong \CO[-1]$, we have
$$
K''' \to (p_Y\times p_X)_*K(e) \exto{\rho'''}
E^*(-y)\boxtimes\CO[2].
$$
Note, that the map $\rho'''$ is related to the map $\rho''$ above
by the following functorial isomorphism
\begin{multline*}
\Hom^\cdot(-,E^*(-y)\boxtimes\CO(x-2e)[3]) =
\Hom^\cdot(-,(\id\times p_X)^!\big(E^*(-y)\boxtimes\CO[2]\big)) = \\ =
\Hom^\cdot((\id\times p_X)_*(-),E^*(-y)\boxtimes\CO[2]),
\end{multline*}
where $(\id\times p_X)^!(F) = \id\times p_X)^*(F)\otimes\CO(x-2e)[1]$
is the right adjoint functor to $(\id\times p_X)_*$.
Therefore $\rho'''\ne 0$.
Note that $(p_Y\times p_X)\cdot i = j$, and $K(e)\cong i_*\CO_W(e)$.
Hence we have the following exact triangle
$$
K''' \to j_*\CO_W(e) \exto{\rho'''} E^*(-y)\boxtimes\CO[2].
$$
with $\rho'''\ne0$, and the functor
$F \mapsto {p_X}_*\Phi_K(R_{\CO(2e-y)}F(e))$ is isomorphic
to the kernel functor $\Phi_{K'''}$. Thus it remains to show
that $\Phi_{K'''}(\CA_Y)=0$.

Using resolution~(\ref{jw}) we deduce that
$\Hom(j_*\CO_W(e),E^*(-y)\boxtimes\CO[2]) \cong
\Hom(E^*(-1),E^*(-1)) = \C$ because $E$ is stable by
proposition~\ref{pm0}. Hence $\rho'''$ comes from the identity
morphism $E^*(-1) \to E^*(-1)$, and $K'''$ is quasiisomorphic to
the complex
$
\CO(-y) \boxtimes V/\CU \to \CO\boxtimes \CU^*.
$
It remains to note that
$
H^\cdot(Y,F\otimes\CO(-y)) = \Hom^\cdot(\CO(y),F) = 0,\quad
H^\cdot(Y,F\otimes\CO) = \Hom^\cdot(\CO,F) = 0,
$
for any $F\in\CA_Y$ by~(\ref{sodec}), hence
$\Phi_{\CO(-y) \boxtimes V/\CU}(\CA_Y) = \Phi_{\CO\boxtimes \CU^*}(\CA_Y)=0$,
hence $\Phi_{K'''}(\CA_Y)=0$.
\end{proof}

\section{Some properties of the category $\CA_Y$}

\subsection*{Serre functor}
Take arbitrary $n,d\in\ZZ$ such that $n+2>d$. 
Let for a moment $Y$ be a smooth $n$-dimensional hypersurface of degree $d$
in $\PP^{n+1}$. Then $Y$ is a Fano manifold and it is easy to check 
that $(\CO_Y,\dots,\CO_Y(n+1-d))$ is an exceptional 
collection in $\D^b(Y)$. Consider the category 
$\CA_Y = \lan\CO_Y,\dots,\CO_Y(n+1-d)\ran^\perp\subset\D^b(Y)$, so that
$$
\D^b(Y) = \lan \CA_Y,\CO_Y,\dots,\CO_Y(n+1-d)\ran
$$
is a semiorthogonal decomposition.

Consider the functor $\SO:\D^b(Y) \to \D^b(Y)$ defined as follows:
$$
\SO(F) = L_\CO (F\otimes\CO_Y(1))[-1].
$$
Note that $\SO$ takes $\CA_Y$ to $\CA_Y$.

\begin{lemma}
We have an isomorhism of functors 
$\SO_{|\CA_Y}^{n+2-d} \cong \FS_{\CA_Y}^{-1}[d-2]$.
\end{lemma}
\begin{proof}
Let $\Phi:\D^b(Y)\to\D^b(Y)$ denote the functor $F\mapsto F(1)$. 
Then using lemma~\ref{philr}, isomorphism 
$\FS_{\D^b(Y)} = \Phi^{d-2-n}[n]$ and proposition~\ref{mut_funct}
we get
\begin{multline*}
\SO_{|\CA_Y}^{n+2-d} = 
(L_{\CO_Y} \circ \Phi[-1]) \circ (L_{\CO_Y} \circ \Phi[-1]) \circ \dots 
\circ (L_{\CO_Y} \circ \Phi[-1])
\cong \\ \cong
L_{\CO_Y} \circ L_{\CO_Y(1)} \circ \dots \circ 
L_{\CO_Y(n-2)} \circ \Phi^{n+2-d}[d-2-n] 
\cong \\ \cong
L_{\lan\CO_Y,\dots,\CO_Y(n-2)\ran} \circ \FS_{\D^b(Y)}^{-1}[d-2] \cong
\FS_{\CA_Y}^{-1}[d-2].
\end{multline*}
\end{proof}

\begin{lemma}
We have an isomorhism of functors $\SO_{|\CA_Y}^d \cong [2-d]$.
\end{lemma}
\begin{proof}
Note that $\SO = \Phi_{K_1}$ with the kernel $K_1$ 
represented by the following complex
$$
\CO_Y(1) \boxtimes \CO \to {\Delta_Y}_* \CO_Y(1).
$$
Iterating, we find that $\SO^d = \Phi_{K_d}$ with the kernel $K_d$ 
represented by the following complex
$$
\CO_Y(1) \boxtimes {\Omega^{d-1}(d-1)}_{|Y} \to \dots \to
\CO_Y(d-1) \boxtimes {\Omega^1(1)}_{|Y} \to
\CO_Y(d) \boxtimes \CO_Y \to 
{\Delta_Y}_* \CO_Y(d).
$$
On the other hand, restricting a resolution of the diagonal 
in $\PP^{n+1}\times\PP^{n+1}$ to $Y\times Y$ we see that
the complex
\begin{multline*}
0 \to
\CO_Y(d-1-n)\boxtimes {\Omega^{n+1}(n+1)}_{|Y} \to \dots \to
\CO_Y    \boxtimes {\Omega^d(d)}_{|Y} \to \\ \to
\CO_Y(1) \boxtimes {\Omega^{d-1}(d-1)}_{|Y} \to \dots \to
\CO_Y(d-1) \boxtimes {\Omega^1(1)}_{|Y} \to
\CO_Y(d) \boxtimes \CO_Y \to 
{\Delta_Y}_* \CO_Y(d)
\end{multline*}
is quasiisomorphic to 
$$
L_1(\alpha\times\alpha)^*\Delta_*\CO_{\PP^{n+1}}(d) \cong {\Delta_Y}_*\CO_Y,
$$
where $\alpha:Y\to\PP^{n+1}$ is the embedding.
Applying the natural morphism between these two complexes we deduce
that ${\Delta_Y}_*\CO_Y[2-d]$ is quasiisomorphic to the complex
$$
\CO_Y(d-1-n)\boxtimes {\Omega^{n+1}(n+1)}_{|Y} \to \dots \to
\CO_Y    \boxtimes {\Omega^d(d)}_{|Y} \to K_d.
$$
It remains to note that
$$
\Phi_{\CO_Y(d-1-n)\boxtimes {\Omega^{n+1}(n+1)}_{|Y}}(\CA_Y) = \dots =
\Phi_{\CO_Y\boxtimes {\Omega^d(d)}_{|Y}}(\CA_Y) = 0,
$$
hence 
${\Phi_{K_d}}_{|\CA_Y} \cong 
{\Phi_{{\Delta_Y}_*\CO[2-d]}}_{|\CA_Y} 
\cong [2-d]$.
\end{proof}

\begin{corollary}
If $c$ is the greatest common divisor of $d$ and $n+2$, then
$\FS_{\CA_Y}^{d/c} \cong [(d-2)(n+2)/c]$.
\end{corollary}

\begin{corollary}
If $Y$ is a cubic threefold then $\FS_{\CA_Y}^3 \cong [5]$.
If $Y$ is a cubic fourfold then $\FS_{\CA_Y} \cong [2]$.
\end{corollary}

%\begin{corollary}
%If $Y$ is a cubic threefold then $\CA_Y$ is not equivalent
%to $\D^b(M)$ for any $M$.
%\end{corollary}

\subsection*{Objects}

Let $Y$ be a smooth cubic threefold.
Simplest examples of objects in $\CA_Y$ are provided by instantons.
\begin{lemma}
If $\CE$ is an instanton of charge~$2$ on $Y$ then
$\CE\in\CA_Y$ and $\CE(-1)\in\CA_Y$.
\end{lemma}
\begin{proof}
Follows from proposition~\ref{charge2}.
\end{proof}

Another examples of objects in $\CA_Y$ 
are provided by curves with theta-characteristics.

\begin{lemma}
Let $M$ be a smooth curve and let $\CL$ be a nondegenerate
theta-characteristic on $M$. For any map $\mu:M\to Y$ the natural morphism 
$H^0(M,\CL\otimes\mu^*\CO(1))\otimes\CO_Y \to (\mu_*\CL)\otimes\CO(1)$
is surjective and its kernel $\CF_{\mu,\CL}\in\CA_Y$.
\end{lemma}
\begin{proof}
Evident.
\end{proof}

Taking $\PP^1$ as a curve, $\CO_{\PP^1}(-1)$ as a theta-characteristic,
and considering omly maps $\mu$ of degree~$1$, we obtain a family of objects 
$\CF_L$ in $\CA_Y$, parameterized by the Fano surface of lines $L$ on $Y$
(in fact, $\CF_L$ is nothing but the sheaf of ideals of $L\subset Y$).
It's Albanese variety is well known to be isomorphic
to the intermediate jacobian of $Y$. 

So, if one would be able to define a notion of stability in $\CA_Y$
in such a way, that any stable object in $\CA_Y$ numerically equivalent
to some $\CF_L$ would be isomorphic to some $\CF_{L'}$, then the Fano
surface would become a moduli space of stable objects in $\CA_Y$,
and it would be possible to reconstruct the intermediate jacobian 
of $Y$ from $\CA_Y$. Since Torelly theorem holds for cubic threefolds 
(see~\cite{CG,T}) it would prove that $\CA_Y$ and $\CA_{Y'}$
are equivalent if and only if $Y\cong Y'$. It would follow also
that $\CA_X$ and $\CA_{X'}$ are equivalent if and only if $X$ 
and $X'$ are birational.

However, it is quite unclear how such stability notion can be defined.

\section*{Appendix A. The Pfaffian hypersurface of a net of skew-forms}\label{s1}

\setcounter{section}{0}
\refstepcounter{section}
\renewcommand{\thesection}{\Alph{section}}
\renewcommand{\thetheorem}{\Alph{section}.\arabic{theorem}}

Let $A=\C^n$ and $V=\C^{2m}$. An {\em $A$-net of skew-forms on $V$}\/
is a linear embedding $f:A\to\Lambda^2V^*$.
Then $F(a)=\Pf(f(a))$ is a homogeneous polynomial of degree $m$ on $A$.
Let $Y=Y_f$ be the corresponding hypersurface of degree $m$ in $\PP(A)$.
We call $Y$ {\em the Pfaffian hypersurface}\/ of the $A$-net~$f$.

The $A$-net $f$ induces a morphism of coherent sheaves on $\PP(A)$
$$
V\otimes\CO_{\PP(A)}(-1) \exto{f} V^*\otimes\CO_{\PP(A)}.
$$
This map is an isomorphism outside of $Y$.
Let $E=E_f$ denote its cokernel. It is a coherent sheaf
on $\PP(A)$ with support on $Y$. We call $E$ {\em the theta-bundle}\/
of the $A$-net. This terminology is suggested by an analogy with
the role of a theta-characteristic on a degeneration curve
of a net of quadrics~\cite{T1}.

Thus for any $A$-net $f$ we have the following exact sequence
\begin{equation}\label{es_app}
0 \to V\otimes\CO(-1) \to V^*\otimes\CO \to \alpha_*E_f \to 0,
\end{equation}
where $\alpha$ is the closed embedding $Y\to\PP(A)$. Second morphism
in this sequence induces an isomorphism
$$
\gamma_f:V^* = H^0(\PP(A),V^*\otimes\CO_{\PP(A)}) \to
H^0(\PP(A),\alpha_*E_f) = H^0(Y_f,E_f).
$$

\begin{definition}
An $A$-net $f$ is called {\em regular}\/ if $\rank f(a) \ge 2m-2$
%the corank of $f(a)$ is not greater then~$2$ 
for any $0\ne a\in A$.
\end{definition}

\begin{remark}
Dimension calculations imply that a regular $A$-net $f$ may exist
only for $\dim A\le 6$.
\end{remark}

\begin{theorem}\label{pff}
Associating to an $A$-net $f$ the triple $(Y_f,E_f,\gamma_f)$
gives a $\GL(A)\times\GL(V)$-equivari\-ant isomorphism between
\begin{itemize}
\item the subset of $\PP(A^*\otimes\Lambda^2V^*)$ formed by all
regular $A$-nets of skew-forms on $V$, and
\item the set of triples $(Y,E,\gamma)$, where
$Y$ is a hypersurface of degree $m$ in $\PP(A)$,
$E$ is a rank~$2$ locally free sheaf on $Y$,
and $\gamma$ is an isomorphism $V^*\to H^0(Y,E)$,
such that
\begin{equation}\label{chern_and_h1}
\begin{array}{c}
c_1(E) = (m-1)[h],\quad c_2(E) = \frac{(m-1)(2m-1)}6h^2,\smallskip\\
H^\cdot(Y,E(t)) = 0\quad\text{for $-(n-2)\le t\le -1$},
\end{array}
\end{equation}
where $[h]\in H^2(Y,\ZZ)$ is the class of a hyperplane section.
\end{itemize}
Further, the theta-bundle $E_f$ of a regular $A$-net is generated
by global sections, $H^0(Y_f,E_f)=V^*$, and induces an embedding
$\kappa:Y_f\to\Gr(2,V)$.
\end{theorem}

\begin{proof}
First of all, we prove that for the theta-bundle of a regular
$A$-net $f$ conditions~(\ref{chern_and_h1}) are satisfied.
This is done by a straightforward calculations, based on
the exact sequence~(\ref{es_app}). This sequence also
implies that $E$ is generated by global sections. It remains to check
that $\kappa$ is an embedding. Note that $Y$ parameterize degenerate
skew-forms in the $A$-net, and $\kappa$ takes a degenerate skew-form
to its kerenl. If two skew-forms have the same kernel, then
a certain linear combination of these skew-forms has $\rank \le 2m-4$, 
which contradicts the regularity of the $A$-net.

Now, assume that $(Y,E,\gamma)$ is a triple,
satisfying~(\ref{chern_and_h1}). Then
$H^\cdot(\PP(A),\alpha_*E(t)) = H^\cdot(Y,E(t))$,
hence it is zero for $-(n-2)\le t \le -1$. Now, let us compute
$H^\cdot(\PP(A),\alpha_*E(t))$ for $t=0$ and $t=-(n-1)$.
To this end choose a line $\PP^1\cong L\subset\PP(A)$
not lying on~$Y$. Then $L\cap Y$ is a $0$-dimensional subscheme
in $Y$ of length $\deg Y = m$. The line $L$ is cut out in $\PP(A)$
by $(n-2)$ hyperplanes, hence $L\cap Y$ is cut out in $Y$
by $(n-2)$ hyperplanes. Therefore we have the Koszul resolution
$$
0 \to E(-(n-2)) \to E(-(n-3))^{\oplus (n-2)} \to \dots \to
E(-1)^{\oplus (n-2)} \to E \to E_{L\cap Y} \to 0.
$$
It follows from~$(\ref{chern_and_h1})$ that
$$
H^p(\PP(A),\alpha_*E) = H^p(Y,E) = H^p(L\cap Y,E_{L\cap Y}) =
\begin{cases}\C^{2m},& \text{for $p=0$}\\0, & \text{for $p>0$}\end{cases}
$$
Twisting the Koszul resolution by $\CO_{\PP(A)}(-1)$ we see that
%\begin{multline*}
$$
H^p(\PP(A),\alpha_*E(-(n-1))) = H^p(Y,E(-(n-1)) \\ 
= H^{p-(n-2)}(L\cap Y,E_{L\cap Y}) =
\begin{cases}\C^{2m},& \text{for $p=n-2$}\\0, & \text{for $p\ne n-2$}\end{cases}
$$
%\end{multline*}
Let us denote $V' = H^{n-2}(\PP(A),\alpha_*E(-(n-1)))$
and recall that we have fixed an isomorphism 
$\gamma:V^* \cong H^0(\PP(A),\alpha_*E)$.
Summarizing, we see that
\begin{equation}\label{h1}
H^p(\PP(A),\alpha_*E(t)) = \begin{cases}
V^*, & \text{for $p=0$, $t=0$}\\
V'  , & \text{for $p=n-2$, $t=-(n-1)$}\\
0  , & \text{for other $(p,t)$ with $-(n-1)\le t\le 0$}
\end{cases}
\end{equation}
Now we can describe the sheaf $\alpha_*E$ via the Beilinson spectral
sequence on $\PP(A)$. It follows from~$(\ref{h1})$ that the spectral
sequence degenerates in the $(n+1)$-th term and gives
$$
0 \to V'\otimes\CO_{\PP(A)}(-1) \exto{f} 
V^*\otimes\CO_{\PP(A)} \to \alpha_*E \to 0.
$$
Dualizing this sequence and twisting it by $\CO_{\PP(A)}(-1)$ we get
$$
0 \to V\otimes\CO_{\PP(A)}(-1) \exto{f^*(-1)} {V'}^*\otimes\CO_{\PP(A)} \to
\CExt^1(\alpha_*E,\CO_{\PP(A)}(-1)) \to 0.
$$
But since $E$ is locally free on $Y$ it follows that
\begin{multline*}
\CExt^1(\alpha_*E,\CO_{\PP(A)}(-1)) \cong 
\alpha_*(E^*\otimes\CHom(L^1\alpha^*\alpha_*\CO_Y,\CO_Y(-1))) \cong \\ \cong
\alpha_*(E^*\otimes\CHom(\CO_Y(-m),\CO_Y(-1))) \cong
\alpha_*(E^*(m-1)).
\end{multline*}
Since $\Lambda^2E = \det E \cong \CO_Y(m-1)$ by $(\ref{chern_and_h1})$
it follows that there exists a skew-symmetric isomorphism
$\sigma:E \to E^*(m-1)$ and $\alpha_*\sigma:\alpha_*E \to \alpha_*E^*(m-1)$.
Since the Beilinson spectral sequence is functorial 
$\alpha_*\sigma$ induces isomorphisms $g:V'\to V$
and $h:V^*\to {V'}^*$ such that 
the following diagram is commutative:
$$
\begin{CD}
0 @>>> V'\otimes\CO_{\PP(A)}(-1) @>{f}>> V^*\otimes\CO_{\PP(A)} 
@>>> \alpha_*E @>>> 0\\
@. @V{g}VV @V{h}VV @V{\alpha_*\sigma}VV \\
0 @>>> V\otimes\CO_{\PP(A)}(-1) @>{f^*}>> {V'}^*\otimes\CO_{\PP(A)} 
@>>> \alpha_*(E^*(m-1)) @>>> 0
\end{CD}
$$
Dualizing this diagram and twisting it by $\CO(-1)$ we get
$$
\begin{CD}
0 @>>> V'\otimes\CO_{\PP(A)}(-1) @>{f}>> V^*\otimes\CO_{\PP(A)} 
@>>> \alpha_*E @>>> 0\\
@. @V{h^*}VV @V{g^*}VV @V{\alpha_*\sigma^*}VV \\
0 @>>> V\otimes\CO_{\PP(A)}(-1) @>{f^*}>> {V'}^*\otimes\CO_{\PP(A)} 
@>>> \alpha_*(E^*(m-1)) @>>> 0
\end{CD}
$$
Now note, that due to skew-symmetry of $\sigma$ we have $\sigma^*=-\sigma$.
It follows that $h^*=-g$ and $g^*=-h$. Identfying
$V'$ with $V$ via $g$ and using the commutativity of
the second diagram we see that
$$
(h f)^*=f^*h^* = g^*f = - h f,
$$
hence
$h f\in\Hom(V\otimes\CO_{\PP(A)}(-1),V^*\otimes\CO_{\PP(A)}) = 
V^*\otimes V^*\otimes A^*$ 
is skew-symmetric with respect to~$V$, therefore is given by an $A$-net
of skew-forms $f:A\to\Lambda^2V^*$. Finally, it is easy to see that
the $A$-net $f$ is regular (because $E$ is locally free of rank~2),
that $Y$ is its Pfaffian hypersurface and that $E$ is its theta-bundle.
\end{proof}

For any $A$-net of skew-forms $f:A\to\Lambda^2V^*$ let $X_f$ denote
the scheme-theoretic intersection of the Grassmannian
$\Gr(2,V)\subset\PP(\Lambda^2V)$ with the codimension $n$
linear subspace $\PP(f(A)^\perp)\subset\PP(\Lambda^2V)$.
%The $A$-net $f$ can be in a unique way reconstructed
%from the subvariety $X_f\subset\Gr(2,V)$.

\begin{proposition}\label{sm}
If $f$ is a regular $A$-net then
$\sing(X_f) = \sing(Y_f) = X_f\cap Y_f \subset \Gr(2,V)$.
In particular, $Y_f$ is smooth iff $X_f$ is smooth.
\end{proposition}
\begin{proof}
Let $U$ be a 2-dimensional subspace of $V$. Then $U$ lies on $X_f$
iff $U$ is isotropic with respect to all skew-forms from the $A$-net $f$.
The tangent space to $\Gr(2,V)$ at $U$ is $\Hom(U,V/U)$. The normal space 
of $\PP(f(A)^\perp)$ in $\PP(\Lambda^2V)$ at $\Lambda^2U$ is $A^*$.
The map $\Hom(U,V/U) \to A^*$ is dual to the map
$A\otimes U \to (V/U)^*$, $(a,u)\mapsto f(a)(u,-)$.
Therefore, $U$ is a singular point of $X_f$ iff $U$ lies in
the kernel of some skew-form from the $A$-net.
Thus $\sing(X_f) = X_f\cap Y_f$.

On the other hand, let $a\in\PP(A)$. Then $a$ lies on $Y_f$ iff $f(a)$
is a degenerate skew-form. Since $f$ is regular, $\rank f(a) = 2m-2$, 
hence its kernel $U$ is 2-dimensional. 
The tangent space to $\PP(A)$ at $a$ is $A/\C a$. The normal space 
of the locus of degenerate skew-forms in $\PP(\Lambda^2V)$ at $f(a)$
is $\Lambda^2U^*$. The map $A/\C a \to \Lambda^2U^*$ is given by
$a'\mapsto f(a')_{|U}$. Therefore, $a$ is a singular point of $Y_f$ 
iff all skew-forms forms from the $A$-net $f$ vanish on $U$. 
Thus, $\sing(Y_f) = X_f\cap Y_f$.
\end{proof}

\section*{Appendix B. Instanton bundles on Fano threefolds of index~$2$}
\refstepcounter{section}

Let $Y$ be a smooth Fano threefold of index~2, so that $\omega_Y=\CO_Y(-2)$.
Let $d=-c_1(\omega_Y)^3/8$ be the degree of $Y$.

\begin{definition}
A sheaf $\CE$ on $Y$ is called {\em instanton bundle}\/ if
$\CE$ is locally free of rank~$2$, stable and
$$
c_1(\CE) = 0,\qquad H^1(Y,\CE(-1))=0.
$$
The {\em topological charge}\/ of an instanton $\CE$ is an integer $k$,
such that $c_2(\CE)=k[l]$, where $[l]\in H^4(Y,\ZZ)$ is the class of a line.
\end{definition}

This definition is a straightforward analog of the definition
of (mathematical) instanton vector bundle on $\PP^3$, see~\cite{OSS}.

\begin{lemma}\label{ci}
If $\CE$ is an instanton vector bundle of charge $k$ on
a Fano threefold~$Y$ of index~$2$ then the dimensions
of the cohomology spaces of twists of $\CE$ are given
by the following table:
$$
\begin{array}{|c|c|c|c|c|c|}
\hline
t & -3 & -2 & -1 & 0 & 1 \\
\hline
h^3(\CE(t)) & \le 2d & 0 & 0 & 0 & 0 \\
\hline
h^2(\CE(t)) & \le 2k-4 & k-2 & 0 & 0 & 0 \\
\hline
h^1(\CE(t)) & 0 & 0 & 0 & k-2 & \le 2k-4 \\
\hline
h^0(\CE(t)) & 0 & 0 & 0 & 0 & \le 2d \\
\hline
\end{array}
$$
where $h^p(\CE(t)) = \dim H^p(Y,\CE(t))$,
and $d$ is the degree of $Y$. Moreover,
$$
h^3(\CE(-3)) = h^0(\CE(1)),\quad
h^2(\CE(-3)) = h^1(\CE(1)),\ \text{and}\
h^0(\CE(1)) - h^1(\CE(1)) = 2d - 2k + 4.
$$
\end{lemma}
\begin{proof}
Note, that $h^0(\CE(t))=0$ for $t\le 0$ by stability
and that Serre duality gives
$$
h^p(\CE(t)) = h^{3-p}(\CE(-2-t))
$$
Hence $h^3(\CE(t))=0$ for $t\ge -2$ and $h^2(\CE(-1))=0$.
Choosing a generic codimension~3 plane section
$L$ of $Y$ we get a Koszul resolution
$$
0 \to \CE(-2) \to \CE(-1)^{\oplus 3} \to \CE^{\oplus 3} \to \CE(1)
\to \CE(1)_{|L} \to 0.
$$
Since $L$ is a $0$-dimensional subscheme in $Y$ we have
$$
H^{>0}(Y,\CE(1)_{|L}) = 0, \qquad
H^0(Y,\CE(1)_{|L}) = r(\CE)\cdot\deg Y = 2d.
$$
Hence the hypercohomology spectral sequence of the Koszul resolution
implies $h^1(\CE(-2))=0$. Then $h^2(\CE)=0$ by Serre duality, hence
$h^2(\CE(1))=0$ by the spectral sequence. Further, we have $h^1(\CE(-3))=0$ 
by Serre duality. And again from the spectral sequence we deduce
$h^1(\CE(1)) \le 2h^1(\CE)$. Finally, using the Riemann-Roch we get
$$
h^1(\CE) = -\chi(\CE) = k-2,\quad
h^0(\CE(1)) - h^1(\CE(1)) = \chi(\CE(1)) = 2d - 2k + 4.
$$
and lemma follows.
\end{proof}

\begin{corollary}\label{charge2_app}
The minimal possible charge for instantons on a smooth Fano threefold
of index~$2$ is $2$, and
%Instantons of charge~$1$ on a smooth Fano threefold~$Y$ of index~$2$
%don't exist. 
if $\CE$ is an instanton of charge~$2$ on $Y$ then
$$
H^p(Y,\CE(t)) = \begin{cases}
\C^{2d}, & \text{for $(p,t)=(0,1)$ and $(p,t)=(3,-3)$}\\
0  , & \text{for other $(p,t)$ with $-3\le t\le 1$}
\end{cases}
$$
\end{corollary}

Following the analogy with instanton bundles on $\PP^3$ we introduce
the following.

\begin{definition}
Let $L\subset Y$ be a line. We say that $L$ is {\em jumping line}\/
for an instanton $\CE$ on $Y$ if $\CE_{|L}\cong\CO_L(t)\oplus\CO_L(-t)$
with $t>0$.
\end{definition}

%One could hope, that as in the $\PP^3$ case the set of jumping lines of
%an instanton on $Y$ is a divisor on the Fano surface of lines on $Y$.

\end{document}